\documentclass{article}
\usepackage{amsmath}
\usepackage{amssymb}
\usepackage{amscd}
\usepackage{amstext}
\usepackage{amsfonts}
\usepackage{euscript}
\usepackage{graphicx}
\usepackage[debug]{hyperref}

\def\scr{\EuScript}

\newcommand{\pcirc}{{\scriptstyle \,\circ\,}}

\newcommand{\C}{\mathbb{C}}
\newcommand{\ZZ}{\mathbb{Z}} \newcommand{\NN}{\mathbb{N}}
\newcommand{\Q}{\mathbb{Q}}

\newcommand{\R}{R}

\newcommand{\W}{\mathbb{W}}
\DeclareMathOperator{\Hom}{\it Hom}

\DeclareMathOperator{\End}{\it End}

\DeclareMathOperator{\ann}{\rm ann}

\DeclareMathOperator{\Dual}{\Bbb D}

\DeclareMathOperator{\Der}{\it Der}

\newcommand{\ux}{\underline{x}}

\newcommand{\uxi}{\underline{\xi}}

\newcommand{\upartial}{\underline{\partial}}

\DeclareMathOperator{\Sim}{{\rm Sym}}
\DeclareMathOperator{\U}{{\rm
U}}

\DeclareMathOperator{\Rees}{{\cal R}}

\DeclareMathOperator{\jac}{\rm Jac}

\newcommand{\Thetafs}{\Theta_{f,s}}

\newcommand{\bul}{\null}

\newcommand{\hol}{{\scr O}}

\newcommand{\derlogD}{\Der(\log D)}
\newcommand{\VO}{{\cal V}_0}

\newcommand{\D}{{\scr D}}
\newcommand{\K}{{\scr K}}
\newcommand{\E}{{\scr E}}

\newcommand{\LL}{{\scr L}}

\newcommand{\M}{{\scr M}}

\newcommand{\OO}{{\scr O}}

\newcommand{\gota}{\mathfrak{a}}

\newcommand{\calV}{{\cal V}}

\newcommand{\Lotimes}{\stackrel{L}{\otimes}}
\DeclareMathOperator{\DR}{DR}

\newcommand{\cd}[3]{{D}^{#1}_{#2}(#3)}
\DeclareMathOperator{\Gr}{Gr}

\newcommand{\OX}{{\scr O}_X}
\newcommand{\DX}{{\scr D}_X}

\newcommand{\DXp}{{\scr D}_{X,p}}
\newcommand{\OXp}{{\scr O}_{X,p}}

\DeclareMathOperator{\SP}{Sp}

\newcommand{\gdltp}{(GDL)$_p$}

\newcommand{\gdlt}{(GDL)}

\newcommand{\gclt}{(GCL)}

\newcommand{\sKp}{(GK)$_p$}

\newcommand{\sK}{(GK)}

\newcounter{numero}[subsection]

\newcounter{snumero}[section]

\renewcommand{\thenumero}{(\thesubsection .\arabic{numero})}
\renewcommand{\thesnumero}{(\thesection .\arabic{snumero})}
\newenvironment{corolario}{\medskip
\refstepcounter{numero}\noindent {\bf  \thenumero\ Corollary.}\
\it}{\vspace{1ex}\par}

\newenvironment{teorema}{\medskip
\refstepcounter{numero}\noindent {\bf \thenumero\ Theorem.}\
\it}{\vspace{1ex}\par}

\newenvironment{definicion}{\medskip
\refstepcounter{numero}\noindent {\bf \thenumero\ Definition.}\
}{\vspace{1ex}\par}

\newenvironment{proposicion}{\medskip
\refstepcounter{numero}\noindent {\bf \thenumero\ Proposition.}\
\it}{\vspace{1ex}\par}

\newenvironment{nota}{\medskip
\refstepcounter{numero}\noindent {\bf \thenumero\ Remark.}\
}{\vspace{1ex}\par}

\newenvironment{cuestion}{\medskip
\refstepcounter{numero}\noindent {\bf \thenumero\ Question.}\ }{\vspace{1ex}\par}

\newenvironment{ejemplo}{\medskip
\refstepcounter{numero}\noindent {\bf \thenumero\ Example.}\
}{\vspace{1ex}\par}

\newcommand\numero{\medskip\refstepcounter{numero}\noindent{\bf\thenumero}\hspace{1em}}

\newenvironment{prueba}{
\noindent {\bf Proof.}\ }{\hfill $\Box$\vspace{1ex}\par}

\title{Linearity conditions on the Jacobian ideal and logarithmic--meromorphic comparison for free divisors}

\author{L. Narv\'{a}ez Macarro\thanks{
Partially supported by MTM2004-07203-C02-01, MTM2007-66929
 and FEDER.}}

\date{}

\begin{document}

\maketitle

\mbox{}\hfill {\em Dedicated to L\^e D\~ung Tr\'ang on his sixtieth birthday}

\begin{abstract}
In this paper we survey the role of D-module theory in the comparison between logarithmic and meromorphic de Rham complexes of integrable
logarithmic connections with respect to free divisors, and we present some new linearity conditions on the Jacobian ideal which arise in this
setting.
\medskip

\end{abstract}

\section*{Introduction}

The comparison between the meromorphic and the logarithmic de Rham complexes of a vector bundle endowed with a logarithmic integrable connection
originally appeared in the case of normal crossing divisors \cite{del_70}.

In \cite{ksaito_log}, the notion of free divisor was introduced, and it produced reasonable logarithmic de Rham complexes with applications to
Singularity theory.

In \cite{cas_mond_nar_96}, a logarithmic-meromorphic comparison (for the trivial bundle) was proved for locally quasi-homogeneous free divisors.
The proof is topological and is based on Grothendieck's comparison theorem.

In the appendix A of \cite{es_vi_86}, for the normal crossing case,
and in \cite{calde_ens}, for general free divisors, D-module theory
enters the scene and the bases for an algebraic treatment of the
logarithmic-meromorphic comparison (at least for free divisors) were
established.

D-module theory has been used for characterizing the logarithmic-me\-ro\-mor\-phic comparison for the trivial bundle with respect to a free
divisor in \cite{cas_ucha_exper}, \cite{torre-45-bis}, \cite{calde_nar_fourier}, and in \cite{calde_nar_fourier}, \cite{calde_nar_lct_ilc} for
general integrable connections, and also for proving algebraically the logarithmic-meromorphic comparison under the additional purely algebraic
hypothesis that the Jacobian ideal is of linear type. This hypothesis had been previously proved for locally quasi-homogeneous free divisors in
\cite{calde_nar_compo}.

In this paper we survey the above results and we highlight some (we believe) new linearity type conditions which appear in many examples and
which deserve further study. These conditions seem to play a role in understanding algebraically the logarithmic-meromorphic comparison, and
also have interesting interactions with Bernstein polynomials and Bernstein functional equations.

Let us now comment on the content of this paper.

Section 1 is a survey on the use of D-module theory in the
logarithmic-meromorphic comparison and on the previous results.

In section 2 we first recall the notions of ``commutative linear type" and ``differential linear type". We also give an algebraic criterion for
checking the (pre)Spencer property. Second, we introduce some generalized linear type properties, which are illustrated in the examples in the
third section, and we study their relationship with ``classical" properties. We end the section with a list of open questions.

In section 3 we have collected several examples. As the reader can guess, the computational complexity is often hard and we comment on the
strategies we have followed.

We thank F.J. Calder\'on Moreno, F.J. Castro Jim\'enez and T. Torelli for many discussions on this subject. We also thank J. Mart\'{\i}n Morales
for computational help. Finally, we thank the referee for his careful reading of our paper and for his comments.

\section{A survey on the logarithmic comparison pro\-blem for free divisors}

\subsection{Notations and basic notions}

Let $X$ be a $n$-dimensional complex analytic manifold and $D\subset
X$ a hypersurface (= divisor), and let us denote by $j: U=X-D
\hookrightarrow X$ the corresponding open inclusion. We denote by
 $\OX$ the sheaf of
holomorphic functions on $X$, ${\cal I}_D \subset \OX$ the ideal of $D$, $\DX$ the sheaf of linear differential operators on $X$ (with
holomorphic coefficients), $\Gr \DX$ the graded ring associated with the filtration $F^{\bul}$ by the order and $\sigma(P)$ the principal symbol
of a differential operator $P$. If $J\subset \DX$ is a left ideal, we denote by $\sigma (J)$ the corresponding graded ideal of $\Gr \DX$. We
denote by $\omega_X$ the sheaf of $n$-differential forms on $X$, which carries a canonical structure of right $\DX$-module. Given a complex of
left $\DX$-modules $\M$, we will denote by $\DR \M$ its de Rham complex:
$$ \DR \M = \R \Hom_{\DX}(\OX,\M) = \left(\omega_X \Lotimes_{\DX} \M\right)[-n].$$
If $\M$ is a single left $\DX$-module, the above definition coincides with the classical one
$$ \DR \M = \M \to \Omega^1_X \otimes_{\OX} \M \to \cdots \to \Omega^n_X \otimes_{\OX} \M.$$

Let us denote by $\jac(D)\subset \OX$ the Jacobian ideal of
$D\subset X$, i.e. the coherent ideal  of $\OX$ whose stalk at any
$p\in X$ is the ideal  generated by $h,\frac{\partial h}{\partial
x_1},\dots,\frac{\partial h}{\partial x_n}$, where $h\in {\scr
O}_{X,p}$ is any reduced local equation of $D$ at $p$ and
$x_1,\dots,x_n\in{\scr O}_{X,p}$ is a system of local coordinates
centered at $p$\footnote{This definition of Jacobian ideal is
independent of the (reduced) equation $h$ and does not match the
usual definition in Singularity Theory, where the Jacobian ideal is
generated by $\frac{\partial h}{\partial x_1},\dots,\frac{\partial
h}{\partial x_n}$ and depends on $h$.}.

For any bounded complex $\K$ of sheaves of $\C$-vector spaces on
$X$, let us denote by $\K^{\vee}= R \Hom_{\C_X}(\K,\C_X)$ its
Verdier dual.

If $A$ is a commutative ring (or a sheaf of commutative rings) and
$M$ an $A$-module, we will denote by $\Sim_A(M)$ its symmetric
algebra. If $I\subset A$ is an ideal, we will denote by $\Rees (I)=
\oplus_{d=0}^\infty I^d t^d \subset A[t]$ its Rees algebra.

The sheaf of meromorphic functions along $D$ is denoted by
$\OX(\star D)$. It is filtered by the invertible sheaves $\OX(rD)$,
$r\geq 0$, formed by meromorphic functions with poles along $D$ of
order at most $r$. As usual, let us also denote by $\OX(-rD)$ de
ideal of $\OX$ of functions vanishing at $D$ with order al least
$r$, $r\geq 0$. We have $\OX(0D)=\OX$, $\OX(-D) = {\cal I}_D$,
$\OX(rD)\otimes_{\OX} \OX(sD) = \OX((r+s)D)$ and
$\OX(-rD)=\OX(rD)^*$ for $r,s\in\ZZ$. A key result in $D$-module
theory is that $\OX(\star D)$ is a left {\em holonomic}
$\DX$-module.

Let us denote by $\Omega^t_X(\star D)$ the sheaf of meromorphic
$t$-forms with poles along $D$, for $t=0,\dots,n$, and
$\Omega^{\bullet}_X(\star D)$ the meromorphic de Rham complex (along
$D$).

We say that a meromorphic $t$-form (with poles along $D$) $\omega$
is logarithmic if both $\omega$ and $d \omega$ have simple poles
along $D$. The sheaf of logarithmic $t$-forms along $D$ is denoted
by $\Omega^t_X(\log D)$. By definition, logarithmic forms along $D$
endowed with the exterior differential give rise to a subcomplex of
$\Omega^{\bullet}_X(\star D)$, which is called {\em logarithmic de
Rham complex (along $D$)} and is denoted by $\Omega^{\bullet}_X(\log
D)$.

The $\OX$-module of (holomorphic) vector fields on $X$, or equivalently, the $\C$-derivations of $\OX$, will be denoted by $\Der_{\C}(\OX)$. A
(holomorphic) vector field $\delta$ on $X$ is called {\em logarithmic (along $D$)} if the ideal of $D$ is fixed by $\delta$, i.e. if it is
tangent to $D$ at the smooth locus. The sheaf of logarithmic vector fields (along $D$) is denoted by $\Der(\log D)$. The bracket of two
logarithmic vector fields is still logarithmic.

Let $x_1,\dots,x_n$ be a system of local coordinates on a open
neighborhood $U$ of a point $p\in D$ and let $h$ be a reduced local
equation of $D$ on $U$. The module of syzygies $(a_0,a_1,\dots,a_n)$
of $(h,h'_{x_1},\dots,h'_{x_n})$ is isomorphic to the sheaf of
logarithmic vector fields on $U$ through the correspondences
\begin{eqnarray*}
&\displaystyle (a_0,a_1,\dots,a_n) \dashrightarrow  \sum_{i=1}^n a_i
\frac{\partial }{\partial x_i},&\\
& \left( -\frac{\delta(h)}{h},\delta(x_1),\dots,\delta(x_n)\right) \dashleftarrow  \delta&
\end{eqnarray*}
and so $\Der(\log D)$ is a coherent module.

The above construction can be performed without choosing local coordinates by considering the map
\begin{equation} \label{eq:global-jac} P \in F^1 \DX \mapsto P(h)
\in \jac(D)
\end{equation}
provided that $D$ has a reduced global equation $h:X\to \C$.

The usual pairing of vector fields and differential forms establishes a perfect $\OX$-duality between $\Der(\log D)$ and $\Omega^1_X(\log D)$
 and then both modules are reflexive \cite{ksaito_log}, cor. (1.7).

We say that $D$ is a {\em free divisor} \cite{ksaito_log} if the $\OX$-module $\Der(\log D)$ is locally free (necessarily of rank $n$), or
equivalently if the $\OX$-module $\Omega^1_X(\log D)$ is locally free (of rank $n$).

Normal crossing divisors, plane curves, free hyperplane arrangements
(e.g. the union of reflecting hyperplanes of a complex reflection
group), discriminants of stable mappings or bifurcation sets of
holomorphic functions are examples of free divisors. More recently,
another interesting source of examples of free divisors has been
studied in \cite{gran-mond-nie-schul}, \cite{buch-mond-06}.

Is $D$ is a free divisor, then for any $t\geq 1$ we have $\Omega^t_X(\log D)= \bigwedge^t \Omega^1_X(\log D)$ \cite{ksaito_log}.

The logarithmic-meromorphic comparison problem for $D$ consists of whe\-ther the inclusion $ \Omega^{\bullet}_X(\log D) \hookrightarrow
\Omega^{\bullet}_X(\star D)$ is a quasi-isomorphism or not, or equivalently by Grothendieck's comparison theorem \cite{gro_rham}, whe\-ther the
canonical morphism $\Omega^{\bullet}_X(\log D) \to \R j_* \C_U$ is an isomorphism (in the derived category) or not.

The meromorphic de Rham complex can be described in terms of
D-module theory as $\Omega^{\bullet}_X(\star D) = \DR \OX(\star D)$.
In sections \ref{subsec:dual} and \ref{subsec:l-m-c} we will explain
how, at least in the case of free divisors, the logarithmic de Rham
complex can be also described in terms of D-module theory (see
corollary \ref{cor:fourier-1}) and how the logarithmic-meromorphic
comparison problem can be translated into a comparison between
certain $\DX$-modules  (see theorem \ref{teo:crit-LCT}).

\subsection{Lie algebroids}

The module $\Der_{\C}(\OX)$ of vector fields on $X$ is at the same
time an $\OX$-module and a sheaf of $\C$-Lie algebras. The interplay
between both structures can be summarized by the equality
$$ [\delta,f\delta'] =
f[\delta,\delta'] + \delta(f)\delta'$$ for any vector fields
$\delta,\delta'$ and any holomorphic function $f$.

In order to understand how logarithmic differential operators (see
section \ref{subsec:ldo}) and logarithmic vector fields are related,
at least in the case of free divisors, and the formal analogies with
vector fields and differential operators, it is useful to consider
the notion of Lie algebroid. A general reference for this notion is
\cite{macken-05}. See also \cite{MR2000m:32015}.

\begin{definicion} A Lie algebroid (on $X$) is an
$\OX$-module $\LL$ endowed with a structure of sheaf of $\C$-Lie
algebras and a $\OX$-linear morphism $\rho:\LL \to \Der_{\C}(\OX)$,
called {\em anchor morphism}, which is also a morphism of sheaves of
Lie algebras and satisfies $ [\lambda,f\lambda'] =
f[\lambda,\lambda'] + \rho(\lambda)(f)\lambda'$ for any local
sections $\lambda,\lambda'$ of $\LL$ and any holomorphic function
$f$.
\end{definicion}

The notion of morphism of Lie algebroids is clear and is left up to the reader.

\begin{ejemplo} \label{ex:Lie-algebroids} 1) The first example of Lie algebroid is $\LL =
\Der_{\C}(\OX)$ with the identity as anchor morphism.\\
2) The sheaf of differential operators of degree $\leq 1$, $F^1 \DX
= \OX \oplus \Der_{\C}(\OX)$, with the projection $F^1 \DX \to
\Der_{\C}(\OX)$ as anchor morphism, is a Lie algebroid.\\
3) Any submodule $\LL \subset \Der_{\C}(\OX)$ which is closed for
the bracket is a Lie algebroid with the inclusion as anchor
morphism. This applies in particular to $\LL = \Der(\log D)$.
\end{ejemplo}

An $\OX$-ring is a sheaf of rings $\scr R$ on $X$ endowed with a
morphism of sheaves of rings $\OX\to {\scr R}$. We say that the
$\OX$-ring $\scr R$ is central over $\C$ if the map $\OX\to {\scr
R}$ sends the constant sheaf $\C_X$ into the center of $\scr R$. For
instance, the sheaf $\DX$ of linear differential operators is an
$\OX$-ring central over $\C$.
\medskip

\begin{definicion} \label{def:admissible}
If $(\LL,\rho)$ is a Lie algebroid on $X$ and $\scr R$ is an
$\OX$-ring central over $\C$, we say that a $\C$-linear map
$\varphi: \LL \to {\scr R}$ is {\em admissible} if the following
properties hold:
\begin{enumerate}
\item[(a)] $\varphi$ is left $\OX$-linear.
\item[(b)] $\varphi([\lambda,\lambda']) =
\varphi(\lambda)\varphi(\lambda')-
\varphi(\lambda')\varphi(\lambda)$ for any local sections
$\lambda,\lambda'$ of $\LL$.
\item[(c)] $\varphi(\lambda) f = f \varphi(\lambda) +
\rho(\lambda)(f) 1_{\scr R}$ for any local section $\lambda$ of
$\LL$ and any holomorphic function $f$.
\end{enumerate}
\end{definicion}

Any Lie algebroid $\LL$ has a functorially associated enveloping
algebra $\U \LL$ (cf. \cite{rine-63} and \cite{calde_nar_ferrara},
th. (1.6)), which is an $\OX$-ring central over $\C$, with a
universal admissible map $\LL \to \U\LL$. The enveloping algebra
$\U\LL$ is filtered with $F^0 \U\LL = \OX$, $F^1 \U\LL = \OX \oplus
\LL$ and its graded algebra is commutative. We have then a canonical
morphism of commutative $\OX$-algebras
\begin{equation} \label{eq:PBW}
 \Sim_{\OX} \LL \to \Gr \U\LL.
\end{equation}

\begin{ejemplo} \label{ex:U-Lie-algebroids} 1) In the case of $\LL = \Der_{\C}(\OX)$ the enveloping algebra is nothing
but the ring $\DX$ of differential operators with the
filtration by the order.\\
2) In the case of $\LL = F^1 \DX$ the enveloping algebra is the
polynomial ring $\DX[s]$, where $s$ is a central variable, endowed
with the total order filtration $F_T^{\bul}$ given by
$$ F_T^k \DX[s] = \sum_{i=0}^k (F^i \DX) s^{k-i},\quad \forall k\geq 0.$$
The universal admissible morphism $F^1 \DX= \OX \oplus
\Der_{\C}(\OX) \to \DX[s]$ is the identity on $\Der_{\C}(\OX)$ and
sends any holomorphic function $f$ to $fs$.\\
3) The case $\LL= \Der(\log D)$ will be treated in the following
section when $D$ is a free divisor.
\end{ejemplo}

We have the following version of the Poincar\'{e}-Birkhoff-Witt
theorem. Its proof can be easily deduced from the global case of
Lie-Rinehart algebras \cite{rine-63}, th. 3.1.

\begin{teorema} \label{teo:PBW} If $\LL$ is a Lie algebroid on $X$ which is locally
free of finite rank as $\OX$-module, then the canonical morphism (\ref{eq:PBW}) is an isomorphism.
\end{teorema}

\numero \label{nume:SP} Let us assume that $\LL$ is a Lie algebroid
on $X$ which is locally free of finite rank as $\OX$-module and let
$\E$ be a left $\U \LL$-module. The {\em
Cartan-Eilenberg-Chevalley-Rinehart-Spencer complex} $\SP_{\LL}\E$
is defined in \cite{calde_nar_lct_ilc}, (1.1.7). It is a complex of
left $\U \LL$-modules and carries an obvious natural augmentation
$$ \varepsilon^0:\SP^{0}_{\LL}\E = \U\LL \otimes_{\OX} \E \to
h^0\left(\SP_{\LL}\E\right) = \E.
$$

One easily sees (\cite{calde_nar_lct_ilc}, prop. (1.1.8)) that if
the left $\U\LL$-module $\E$ is locally free of finite rank over
$\OX$ then $\SP_{\LL}\E$ is a $\U\LL$-locally free resolution of the
$\U\LL$-module $\E$, that will be simply called {\em Spencer
resolution} of $\E$.

\subsection{Logarithmic differential operators} \label{subsec:ldo}

The Malgrange-Kashiwara $V$-filtration with respect to $D$ on the sheaf $\DX$ of linear differential operators on $X$ is defined by\footnote{It
was originally defined only for smooth $D$.}
$$ \calV_k^D \DX = \{ P \in \DX\ |\ P({\cal I}_D^m) \subset {\cal
I}_D^{m-k}\ \forall m \in\ZZ \},\quad k\in \ZZ.$$ It is an
increasing filtration and $\calV_0^D \DX$ is a subsheaf of rings of
$\DX$ with $F^0 \calV_0^D \DX = \OX$ and $F^1 \calV_0^D \DX = \OX
\oplus \Der(\log D)$. We have then a canonical morphism of
commutative $\OX$-algebras
\begin{equation} \label{eq:Gr-V0}
 \Sim_{\OX} \Der(\log D) \to \Gr_F \calV_0^D \DX
\end{equation}
and a canonical morphism of non commutative $\OX$-rings
\begin{equation} \label{eq:U-V0}
 \U \Der(\log D) \to  \calV_0^D \DX.
\end{equation}

The $0$'th term of the Malgrange-Kashiwara filtration with respect
to $D$ on the sheaf $\DX$, $\calV_0^D \DX$, is also denoted by
$\DX(\log D)$, or simply by $\VO$ if no confusion is possible, and
its sections are called logarithmic differential operators with
respect to $D$.

The following theorem is not explicitly stated in \cite{calde_ens},
but it is essentially contained in proposition 2.2.5 of {\em loc.
cit.} (see also \cite{schul-criter}).

\begin{teorema} \label{teo:calde-ens} If $D$ is a free divisor, then the canonical
morphism (\ref{eq:U-V0}) is an isomorphism.
\end{teorema}

The proof of the above theorem given in \cite{calde_ens} is based on
the corollary 2.1.6 of {\em loc. cit.}, which states that if $D$ is
free the morphism (\ref{eq:Gr-V0}) is an isomorphism of commutative
$\OX$-algebras. Let us note that, for proving corollary 2.1.6, it is
possible to avoid theorem 2.1.4 and to obtain it directly from
proposition 2.1.2 in \cite{calde_ens}.

As a consequence of theorem \ref{teo:calde-ens} we deduce that if
$D$ is free, then $\VO$ is coherent and it has noetherian stalks of
finite global homological dimension. If
$\{\delta_1,\dots,\delta_n\}$ is a local basis of the logarithmic
vector fields on a connected open set $V$, any differential operator
in $\Gamma(V,\VO)$ can be written locally in a unique way as a
finite sum
\begin{equation} \label{eq:estruc-V0}
\sum_{\substack{\alpha \in\NN^n\\
|\alpha|\leq d}} a_{\alpha} \delta_1^{\alpha_1}\cdots \delta_n^{\alpha_n},\end{equation} where the $a_{\alpha}$ are holomorphic functions on
$V$.

\subsection{Integrable logarithmic connections with respect to a free divisor}

Let $\E$ be a $\OX$-module.  A {\em logarithmic connection} with
respect to $D$ on $\E$ is a $\C$-linear morphism $\nabla:
\E\xrightarrow{} \Omega^1_X(\log D)\otimes_{\OX}\E$ satisfying
Leibniz's rule, i.e. $\nabla(ae) = a \nabla(e) + da\otimes e$ for
any holomorphic function $a$ and any section $e$ of $\E$. The data
of a such logarithmic connection is equivalent to the data of a left
$\OX$-linear map $\nabla_\bullet:\derlogD \xrightarrow{}
\End_{\C}(\E)$, where the Leibniz rule becomes $\nabla_\delta(ae) =
a\nabla_\delta(e) + \delta(a)e$ for any logarithmic vector field
$\delta$, any holomorphic function $a$ and any section $e$ of $\E$.

From now on let us assume that $D$ is a free divisor and that
$\nabla$ is a logarithmic connection on $\E$. As in the classical
case (cf. \cite{del_70}, chap. I, 2.10) we define a sequence of
$\C$-linear morphisms
\begin{equation} \label{eq:pre-log-DR}
\E \xrightarrow{\nabla} \Omega^1_X(\log D)\otimes\E \xrightarrow{\nabla} \Omega^2_X(\log D)\otimes\E \xrightarrow{\nabla} \cdots
\xrightarrow{\nabla} \Omega^n_X(\log D)\otimes\E.
\end{equation}

We say that $\nabla$ is integrable if the sequence (\ref{eq:pre-log-DR}) is a complex of sheaves of complex vector spaces. As in the classical
case cf. \cite[chap.~I, 2.12--2.14]{del_70}, $\nabla$ is integrable if and only if the morphism $\nabla_\bullet$ preserves Lie brackets (i.e. it
is admissible in the sense of definition \ref{def:admissible} ), or merely if $\nabla^2: \E \to \Omega^2_X(\log D)\otimes\E$ vanishes.

If $\nabla$ is integrable then the complex (\ref{eq:pre-log-DR}) will be called the {\em logarithmic de Rham complex} of $\E$ (endowed with
$\nabla$) and will be denoted by $\Omega_X^{\bullet}(\log D)(\E)$.

The following proposition is a straightforward consequence of
theorem \ref{teo:calde-ens} (see \cite{calde_ens}, cor. 2.2.6) and
it is similar to the well known case of integrable connections and
left $\DX$-module structures.

\begin{proposicion} To give an integrable logarithmic connection on
$\E$ is equivalent to giving a structure of left $\DX(\log
D)$-module on $\E$ extending its structure of $\OX$-module.
\end{proposicion}

The following theorem is the first step in describing logarithmic de
Rham complexes in terms of D-module theory (see corollary
\ref{cor:fourier-1}). It is a consequence of theorem
\ref{teo:calde-ens} and of the existence of the Spencer resolution
of $\OX$ \ref{nume:SP}.

\begin{teorema}\label{teo:dif-carac-DRlog}(\cite{calde_ens}, cor.~3.2.2)
For any left $\VO$-module $\E$ there is a ca\-no\-nical isomorphism
in the derived category $ \R \Hom_{\VO}(\OX,\E) \simeq
\Omega_X^{\bullet}(\log D)(\E)$.
\end{teorema}

From now on we will only consider integrable logarithmic connections
on $\OX$-modules which are locally free of finite rank, and such an
$\OX$-module endowed with an integrable logarithmic connection (with
respect to $D$) will be simply called an {\em ILC} (with respect to
$D$). In other words, an ILC will be a left $\VO$-module which is
locally free of finite rank over $\OX$.

The first examples of ILC are the invertible $\OX$-modules $\OX(mD)\subset \OX(\star D)$, $m\in\ZZ$. If $f=0$ is a reduced local equation of $D$
at $p\in D$ and $\delta_1,\dots,\delta_n$ is a local basis of $\Der(\log D)_p$ with $\delta_i(f)=\alpha_i f$, then $f^{-m}$ is a local basis of
$\OO_{X,p}(mD)$ over $\OO_{X,p}$ and from (\ref{eq:estruc-V0}) we have the following local presentation over $\D_{X,p}(\log D)$
\begin{equation*}
\OO_{X,p}(mD) \simeq \D_{X,p}(\log D)/\D_{X,p}(\log D)(\delta_1+m\alpha_1,\dots, \delta_n+m\alpha_n).\end{equation*}

For $\E=\OX$ we have $\Omega_X^{\bullet}(\log D)(\OX)=\Omega_X^{\bullet}(\log D)$.

For any ILC $\E$ and any integer $m$, the locally free $\OX$-modules
$\E(mD):= \E \otimes_{\OX} \OX(mD)$ and $\E^* := \Hom_{\OX}(\E,\OX)$
are endowed with a natural structure of left $\VO$-module (cf.
\cite{calde_nar_fourier}, \S 2 and \cite{calde_nar_ferrara}, \S 2),
and they are again ILC, and the usual isomorphisms $ \E(mD)(m'D)
\simeq \E((m+m')D)$, $\E(mD)^* \simeq \E^*(-mD)$ are $\VO$-linear.

For any ILC $\E$ the complex $\SP\E:=\SP_{\Der(\log D)}\E$ is a
resolution of $\E$ as left $\VO$-module. In particular, any ILC is
coherent as left $\VO$-module.

\subsection{Koszul and Spencer properties}

In this section we recall two natural properties which appear when
studying free divisors and the logarithmic-meromorphic comparison
problem.

From now on $D\subset X$ is assumed to be a free divisor.

\begin{definicion} (\cite{calde_ens}, def. 4.1.1)
The (free) divisor $D$ is said to be  {\em Koszul}
 at a point $p\in D$ if the symbols of any (or
some) local basis $\{\delta_1,\dots,\delta_n\}$ of $\Der(\log D)_p$
form a regular sequence in $\Gr {\scr D}_{X,p}$. We say that $D$ is
a {\em Koszul} divisor if it is so at any point $p\in D$.
\end{definicion}

For the free divisor $D$, to be Koszul is equivalent to being
holonomic in the sense of \cite{ksaito_log}, def. (3.8), i.e. the
logarithmic stratification of $D$ is locally finite\footnote{This
was noticed by M. Schulze.}.

Any plane curve is a Koszul free divisor (see the proof of corollary 4.2.2 in \cite{calde_ens}). Any locally quasi-homogeneous\footnote{i.e. for
which for any $p\in D$ there is a system of local coordinates  $\ux$ centered at $p$ such that the germ $(D,p)$ has a reduced weighted
homogeneous defining equation (with strictly positive weights) with respect to $\ux$.} free divisor $D$ is Koszul \cite{calde_nar_compo}.
Examples of locally quasi-homogeneous free divisors are free hyperplane arrangements or discriminants of stable maps in Mather's ``nice
dimensions".

\begin{definicion}\label{def:spencer} (\cite{cas_ucha_stek}, def. 3.3) The (free) divisor $D$ is said
to be {\em Spencer} if the complex of $\DX$-modules
$\DX\otimes_{\VO} \SP\OX$ is holonomic and concentrated in degree 0.
\end{definicion}

Since $\SP\OX$ is a locally free resolution of the $\derlogD$-module
$\OX$, we have $\DX\Lotimes_{\VO} \OX= \DX\otimes_{\VO} \SP\OX$ and
so $D$ is Spencer if and only if $h^i\left(\DX\Lotimes_{\VO}
\OX\right)=0$ for $i\neq 0$ and $h^0 \left(\DX\Lotimes_{\VO}
\OX\right)=\DX\otimes_{\VO} \OX = \DX/\DX \Der(\log D)$ is
holonomic.
\medskip

\begin{definicion} \label{def:pre-spencer}
The (free) divisor $D$ is said to be {\em pre-Spencer } if the
complex of $\DX$-modules $\DX\otimes_{\VO} \SP\OX$ is concentrated
in degree 0.
\end{definicion}
\medskip

We have the following result (\cite{calde_nar_fourier}, prop. 1.2.3).

\begin{proposicion} Any Koszul free divisor is Spencer.
\end{proposicion}

There are lots of examples of Spencer free divisors which are not Koszul (see for instance the examples in section \ref{sec:ex}).

It should be noted that (pre)Spencer property is not easy to check
algebraically and in principle one needs to check that
\underline{each} $h^i\left(\DX\otimes_{\VO} \SP\OX\right)$ vanishes,
for $i=1,\dots,n$. However, see proposition \ref{prop:crit-sp},
corollary \ref{cor:crit-SP} and corollary \ref{cor:carac-top-SP}.

\subsection{Duality} \label{subsec:dual}

In this section we review some results of \cite{calde_nar_fourier}
on the relationship between the duality of integrable logarithmic
connections with respect to a free divisor $D\subset X$ and the
duality of D-modules.
\medskip

Let us denote by $\Dual_{\DX}: \cd{b}{coh}{\DX}\to \cd{b}{coh}{\DX}$
the duality functor of D-module theory (cf. \cite{meb_formalisme},
def. I.4.1.6). It is defined by
$$\Dual_{\DX} \M = \Hom_{\OX}(\omega_X,\R\Hom_{\DX}(\M,\DX))[n].$$
The duality functor $\Dual_{\DX}$ is a contravariant involutive
self-equivalence of the derived category $\cd{b}{coh}{\DX}$ of
bounded complexes of left $\DX$-modules with coherent homologies,
and it induces a contravariant involutive self-equivalence of the
abelian category of left holonomic $\DX$-modules.
\medskip

If we start with an ILC (with respect to $D$) $\E$, we can take first its dual (as ILC) and second the scalar extension from $\VO$ to $\DX$, and
we obtain $\DX\Lotimes_{\VO} \E^*$, or we can reverse the order and take first its scalar extension and second its dual (as $\DX$-module), and
we obtain $\Dual_{\DX}(\DX\Lotimes_{\VO} \E)$. How are both results related? The following theorem gives a precise answer to this question.

\begin{teorema} (\cite{calde_nar_fourier}, cor. 3.1.2) Let $\E$ be an ILC
(with respect to $D$). There is a natural isomorphism in
$\cd{b}{coh}{\DX}$:
$$ \Dual_{\DX}\left(\DX\Lotimes_{\VO} \E\right)\simeq \DX\Lotimes_{\VO} \E^*(D).$$
\end{teorema}

As a corollary of the above theorem and of theorem
\ref{teo:dif-carac-DRlog} we obtain the following description of
logarithmic de Rham complexes in terms of the de Rham functor in
D-module theory.

\begin{corolario}\label{cor:fourier-1} (\cite{calde_nar_fourier}, cor. 3.1.5) Let $\E$ be an ILC (with respect to $D$).
There is a natural isomorphism in the derived category of bounded
complexes of sheaves of complex vector spaces $\cd{b}{}{\C_X}$:
$$
\Omega^{\bullet}_X(\log D)(\E)  \simeq \DR \left(\DX\Lotimes_{\VO}
\E(D)\right).$$
\end{corolario}

The following theorem is proved in \cite{calde_nar_fourier}, cor.
3.1.6 and 3.1.8. It uses the deep properties of the de Rham functor
(cf. \cite{meb_formalisme}, th. II.4.1.5, th. I.10.13), the local
duality theorem in D-module theory (cf. \cite{meb_formalisme}, th.
I.4.3.1; see also \cite{nar-ldt}) and the faithful flatness of
$\DX^\infty$ over $\DX$ (\cite{skk}; see also \cite{nar-rojas}).
\medskip

\begin{teorema} \label{teo:fourier-2} Let $\E$ be an ILC (with respect to $D$).
The following properties are equivalent:
\begin{enumerate}
\item[1)] The complex $\DX\Lotimes_{\VO} \E(D)$ is holonomic and concentrated in degree $0$.
\item[2)] The complex $\DX\Lotimes_{\VO} \E^*$ is holonomic and concentrated in degree $0$.
\item[3)] The complex $\Omega^{\bullet}_X(\log D)(\E)$ is a perverse sheaf.
\item[4)] The complex $\Omega^{\bullet}_X(\log D)(\E^*(-D))$ is a perverse sheaf.
\end{enumerate}
Moreover, if the above properties hold there is a natural
isomorphism of perverse sheaves $ \Omega_X^{\bullet}(\log D)(\E)
\simeq \Omega_X^{\bullet}(\log D)(\E^*(-D))^{\vee}$.
\end{teorema}

The formula $ \Omega_X^{\bullet}(\log D)(\E) \simeq
\Omega_X^{\bullet}(\log D)(\E^*(-D))^{\vee}$ has been proved in
\cite{es_vi_86}, (A.2) in the case of a normal crossing divisor.
\medskip

The following corollary deals with the case of the trivial ILC
$\E=\OX$.
\medskip

\begin{corolario} \label{cor:carac-top-SP} The following properties are equivalent:
\begin{enumerate}
\item The divisor $D$ is Spencer.
\item The logarithmic de Rham complex $\Omega^{\bullet}_X(\log
D)$ is a perverse sheaf.
\item The complex $\Omega^{\bullet}_X(\log
D)(\OX(-D))$ is a perverse sheaf.
\end{enumerate}
\end{corolario}

\begin{nota} In the above corollary, if $h=0$ is a local reduced equation of $D$, then the complex $\Omega^{\bullet}_X(\log
D)(\OX(-D))$ is nothing but $h \Omega^{\bullet}_X(\log D)$, which has been used in \cite{mond-plms-2000} to describe the Gauss-Manin connection
on the cohomology of families of free divisors.
\end{nota}

\subsection{Logarithmic-meromorphic comparison} \label{subsec:l-m-c}

In this section we assume that $D$ is a free divisor and $\E$ is an ILC along $D$.

Let us denote by $\DX(\star D)$ the sheaf of meromorphic linear differential operators with poles along $D$. One has obvious left and right
$\OX(\star D)$-linear isomorphisms
$$ \OX(\star D)\otimes_{\OX} \DX \stackrel{\text{left}}{\simeq}
\DX(\star D) \stackrel{\text{right}}{\simeq} \DX \otimes_{\OX} \OX(\star D).$$The induced maps $ \OX(\star D)\otimes_{\OX} \VO \xrightarrow{}
\DX(\star D) \xleftarrow{} \VO \otimes_{\OX} \OX(\star D)$ are also isomorphisms and so, if $\E$ is a left $\VO$-module, the localization $
\E(\star D) := \OX(\star D)\otimes_{\OX} \E = \DX(\star D)\otimes_{\VO} \E$ is a left $\DX(\star D)$-module, and by scalar restriction, a left
$\DX$-module. Moreover, if $\E$ is an ILC, then $\E(\star D)$ is a meromorphic connection (locally free of finite rank over $\OX(\star D)$) and
so it is a holonomic $\DX$-module (cf. \cite{meb_nar_dmw}, Th. 4.1.3). Actually, $\E(\star D)$ has regular singularities on the smooth part of
$D$ (it has logarithmic poles! \cite{del_70}) and so it is regular everywhere \cite{meb-cimpa-2}, Cor. 4.3-14, which means that if $\LL$ is the
local system of horizontal sections of $\E$ on $U=X-D$, the canonical morphism $ \DR \E(\star D)= \Omega_X^{\bullet}(\E(\star D)) \xrightarrow{}
R j_* \LL$ is an isomorphism in the derived category.

We say that $\E$ satisfies the logarithmic comparison theorem (LCT)
if the canonical morphism $\Omega_X^{\bullet}(\log D)(\E) \to R j_*
\LL$ is an isomorphism in the derived category. If $\OX$ satisfies
the LCT, we simply say that the divisor $D$ satisfies the LCT.

In \cite{cas_mond_nar_96} it has been proved that the LCT holds for any free divisor $D$ which is locally quasi-homogeneous.

For any ILC $\E$ and any integer $m$, $\E(mD)$ is a sub-$\VO$-module of the regular meromorphic connection (and holonomic $\DX$-module)
$\E(\star D)$, and so we have a canonical morphism in the derived category of left $\DX$-modules
\begin{equation*}
 \rho_{\E,m}: \DX \Lotimes_{\VO} \E(mD)\to \E(\star D),\end{equation*}
given by $\rho_{\E,m}(P\otimes e) = P e$.

We have the following characterization for the LCT (see
\cite{calde_nar_fourier}, th. 41 and \cite{calde_nar_lct_ilc}, th.
(2.1.1)).

\begin{teorema} \label{teo:crit-LCT} Let $\E$ be an ILC (with respect to the free divisor $D$)
and let $\LL$ be the local system of its horizontal sections on $U=X-D$. The following properties are equivalent:
\begin{enumerate}
\item[1)] $\E$ satisfies the LCT.
\item[2)] The inclusion $\Omega_X^{\bullet}(\log D)(\E) \hookrightarrow
\Omega_X^{\bullet}(\E(\star D))$ is a quasi-isomorphism.
\item[3)] The morphism $\rho_{\E,1}: \DX\Lotimes_{\VO} \E(D)\to \E(\star D)$
is an isomorphism in the derived category of left $\DX$-modules.
\item[4)] The complex $\DX\Lotimes_{\VO} \E(D)$ is
concentrated in degree $0$ and the $\DX$-module $\DX\otimes_{\VO} \E(D)$ is holonomic and isomorphic to its localization along $D$.
\item[5)] The canonical morphism $j_!\LL^{\vee} \to
\Omega_X^{\bullet}(\log D)(\E^*(-D))$ is an isomorphism in the derived category of complexes of sheaves of complex vector spaces\footnote{In
other words, $\Omega_X^{\bullet}(\log D)(\E^*(-D))|_D = 0$, i.e. the complex $\Omega_X^{\bullet}(\log D)(\E^*(-D))$ is exact on any point of
$D$.}.
\end{enumerate}
\end{teorema}

In the case of the trivial ILC $\E=\OX$, we have the following.

\begin{corolario} \label{cor:carac-LCT-OX} The following properties are equivalent:
\begin{enumerate}
\item[1)] $D$ satisfies the LCT.
\item[2)] The inclusion $\Omega_X^{\bullet}(\log D) \hookrightarrow
\Omega_X^{\bullet}(\star D)$ is a quasi-isomorphism.
\item[3)] The morphism $\rho=\rho_{\OX,1}: \DX\Lotimes_{\VO} \OX(D)\to \OX(\star D)$
is an isomorphism in the derived category of left $\DX$-modules.
\item[4)] The complex $\DX\Lotimes_{\VO} \OX(D)$ is
concentrated in degree $0$ and the $\DX$-module $\DX\otimes_{\VO} \OX(D)$ is holonomic and isomorphic to its localization along $D$.
\item[5)] The canonical morphism $j_!\C_U \to
\Omega_X^{\bullet}(\log D)(\OX(-D))$ is an isomorphism in the
derived category of complexes of sheaves of complex vector
spaces\footnote{In other words, $\Omega_X^{\bullet}(\log
D)(\OX(-D))|_D=0$, i.e. the complex $\Omega_X^{\bullet}(\log
D)(\OX(-D))$ is exact on any point of $D$.}.
\end{enumerate}
\end{corolario}

\begin{nota} \label{nota:carac-LCT-explicit} Let $f$ be a reduced local equation of $D$ and $\delta_1,\dots,\delta_n$ a local basis of $\derlogD$ with $\delta_i(f)=\alpha_i f$.
Property 3) in corollary \ref{cor:carac-LCT-OX} means the
conjunction of the following properties:
\begin{enumerate}
\item[a)] $D$ is pre-Spencer.
\item[b)] (b-1) The $\DX$-module $\OX(\star D)=\OX(f^{-1})$ is generated by $f^{-1}$ and (b-2) $\ann_{\DX} f^{-1}$ is generated by order one
operators, and so by $\delta_1+\alpha_1,\dots, \delta_n+\alpha_n$.
\end{enumerate}
Actually, after \cite{torre-45-bis}, prop. 1.3, property (b-2) implies property (b-1) and so we conclude that $D$ satisfies the LCT if and only
if it is pre-Spencer and $\ann_{\DX} f^{-1}=\DX(\delta_1+\alpha_1,\dots, \delta_n+\alpha_n)$ (see also \cite{cas_ucha_exper}).
\end{nota}

\begin{nota} \label{nota:LWQH->LCT} In the case of locally quasi-homogeneous free divisors, property 5) in corollary \ref{cor:carac-LCT-OX} holds by \cite{mond-plms-2000},
Lemma 3.3, (6) and so we obtain another proof of the LCT for such divisors \cite{calde_nar_lct_ilc}, cor. (2.1.4). Let us also note that the
same argument can be applied to the case of locally weakly quasi-homogeneous (LWQH) divisors, as it is noted in \cite{cas_gago_hartillo_ucha},
remark 3.11, even without the hypothesis of being Spencer. In fact, the corollary above says in particular that any LWQH free divisor is Spencer
(see also corollary \ref{cor:carac-top-SP}).
\end{nota}

\section{Linearity conditions on the Jacobian ideal}

\subsection{The ring $\D[s]$ and the Bernstein construction} \label{subsec:D-bernstein}

Let us suppose that the divisor $D$ is given by a (locally) reduced global equation $f:X\to \C$.
\medskip

We know from the Poincar\'e-Birkhoff-Witt theorem \ref{teo:PBW} and example \ref{ex:U-Lie-algebroids}, 2), that the canonical map
$\eta:\Sim_{\OX} (F^1\DX) \to \Gr_{F_T} \DX [s]$ is an isomorphism of graded $\OX$-algebras.

The free module of rank one over $\OX[f^{-1},s]$ generated by the symbol $f^s$, $\OX[f^{-1},s]f^s$, has a natural left module structure over
$\DX[s]$: the action of a derivation $\delta\in\Der_{\C}(\OX)$ is given by $\delta (f^s) = \delta (f) s f^{-1} f^s$ (see \cite{bern_72}).

Let us call $\varphi_0:\Sim_{\OX} (F^1\DX) \to \Rees (\jac(D))\subset \OX[t]$ the composition of the canonical surjective map $\Sim_{\OX}
\jac(D) \to \Rees (\jac(D))$ with the surjective map $\Sim_{\OX} (F^1\DX) \to \Sim_{\OX} (\jac(D)) $ induced by (\ref{eq:global-jac}) and
\begin{equation} \label{eq:varphi} \varphi := \varphi_0 \pcirc \eta^{-1}: \Gr_{F_T} \DX
[s] \to \Rees (\jac(D)).\end{equation}

For each $P\in\DX[s]$ of total order $d$, we have that $ P(f^s)=Q(s) f^{-d} f^s$ where $Q(s)$ is a polynomial of degree $d$ in $s$ with
holomorphic coefficients. Let us denote by $C_{P,d}\in \OX$ the leading coefficient of $Q(s)$. A straightforward computation (cf.
\cite{yano_78}, chap.~I, Prop.~2.3) shows that $\varphi (\sigma_T(P))= C_{P,d} t^d$, where
$\sigma_T$ denotes the symbol with respect to the
total order filtration, and so $ \sigma_T(\ann_{\DX[s]} f^s) \subset \ker \varphi$.

It is clear that $F^0_T \ann_{\DX[s]} f^s = 0$ and that
\begin{equation} \label{eq:thetafs}
\Thetafs:= F^1_T \ann_{\DX[s]} f^s
\end{equation}
 is formed by the
operators $\delta - \alpha s$ with $\delta\in \Der_{\C}(\OX)$, $\alpha\in \OX$ and $\delta(f)=\alpha f$. One easily sees that the $\OX$-linear
map
\begin{equation} \label{eq:derlog-thetafs}
\textstyle \delta \in \Der(\log D) \mapsto \delta - \frac{\delta(f)}{f} s\in \Thetafs
\end{equation}
is an isomorphism of Lie algebroids on $X$. We obtain a canonical isomorphism $ \Thetafs \simeq \Gr^1_{F_T} \ann_{\DX[s]} f^s$. Let us denote by
$\ann^{(1)}_{\DX[s]}(f^s)$ the left ideal of $\DX[s]$ generated by the total order one operators in $\ann_{\DX[s]} f^s$, i.e. $
\ann^{(1)}_{\DX[s]}(f^s) =\DX[s]\cdot\Thetafs$. On the other hand, the homogeneous part of degree one $[\ker \varphi]_1\subset\ker \varphi$ is
also canonically isomorphic to $\Thetafs$, and so we obtain $$ \Gr^1_{F_T} \ann_{\DX[s]} f^s \left( = \left[ \sigma_T(\ann_{\DX[s]}
f^s)\right]_1=\sigma_T(\Thetafs) \right) = [\ker \varphi]_1.$$

\subsection{Divisors of linear type}

In this section $D$ will be a hypersurface (= divisor) of $X$.

\begin{definicion}
(Cf. \cite{vascon_cmcaag}, {\S}7.2) Let $A$ be a commutative ring
and $I\subset A$ an ideal. We say that $I$ is of {\em linear type}
if the canonical (surjective) map of graded $A$-algebras $ \Sim_A(I)
\to \Rees(I)$ is an isomorphism.
\end{definicion}
\medskip

Ideals generated by a regular sequence are the first example of ideals of linear type. If $A=\C[x_1,x_2]$ and $f\in A$ is a reduced polynomial
such that the divisor $D= \{f=0\}$ is not locally quasi-homogeneous, then the Jacobian ideal $(h,\frac{\partial h}{\partial x_1},\frac{\partial
h}{\partial x_2})$ is not of linear type (see proposition \ref{prop:carac-2-dim}).
\medskip

\begin{definicion} (\cite{calde_nar_lct_ilc}, def. (1.4.2))
 We say that the divisor $D$ is of {\em commutative
linear type} at $p\in D$ if the stalk at $p$ of its Jacobian ideal is of linear type. We say that $D$ is of {\em commutative linear type} if it
is so at any $p\in D$.
\end{definicion}

\begin{definicion} We say that the divisor $D$ is Euler homogeneous at $p\in D$ if for some (and hence any) reduced local equation $f$ of $D$ at $p$
there is a germ of vector field $\chi$ at $p$ such that $\chi(f)=f$. We say that $D$ is Euler homogeneous if it is so at any $p\in D$.
\end{definicion}

\begin{ejemplo} If $D$ is Euler homogeneous and has isolated singularities, then $D$ is of commutative linear type.
\end{ejemplo}

The following proposition is proved in \cite{calde_nar_lct_ilc},
remark (1.6.6).

\begin{proposicion} \label{prop:clt->E-hom} If $D$ is of commutative linear type (at $p$),
then $D$ is Euler homogeneous (at $p$).
\end{proposicion}

Theorem \ref{teo:crit-LCT} has been used in
\cite{calde_nar_lct_ilc}, \S 3 to give a criterion for the LCT to be
satisfied by integrable logarithmic connections with respect to a
free divisor of commutative linear type.

\begin{nota} \label{nota:expli-clt}
To say that a divisor $D$ given by a reduced global equation $f$ is of commutative linear type is equivalent to saying that $\ker \varphi$ (see
(\ref{eq:varphi})) is generated by its homogeneous part of degree 1, $[\ker \varphi]_1 = \sigma_T(\Thetafs)$ (see (\ref{eq:thetafs})).
\end{nota}
\smallskip

We have the following theorem (\cite{calde_nar_compo}, th.~5.6).

\begin{teorema} \label{teo:lqh->clt} Any locally quasi-homogeneous free divisor is of
commutative linear type.
\end{teorema}

\begin{definicion} Let $p\in D$ and let us write $\hol=\OXp$ and $\D=\DXp$.
We say that $D$ is of {\em differential linear type} at $p\in D$ if for some (and hence any)\footnote{One easily sees that this condition does
not depend on the choice of the local equation.} reduced local equation $f\in\hol$ of $D$ at $p$, the ideal $\ann_{\D[s]} f^s$ is generated by
total order one operators, i.e. $ \ann_{\D[s]} f^s =\ann^{(1)}_{\D[s]} f^s= \D[s]\cdot \Thetafs$ (see (\ref{eq:thetafs})). We say that $D$ is of
{\em differential linear type} if it is so at any $p\in D$.
\end{definicion}

It is clear that being of commutative or differential linear type
for a divisor are open conditions.

The following result is proved in \cite{calde_nar_compo}, prop.~3.2.

\begin{proposicion} \label{prop:clt->dlt} If the divisor $D$ is of commutative
linear type (at $p\in D$), then it is of  differential linear type
(at $p\in D$) and if $f:X\to \C$ is a reduced equation of $D$, then
$ \Gr_{F_T} \ann_{\DX[s]} f^s \left( = \sigma_T(\ann_{\DX[s]}
f^s)\right) = \ker \varphi$.
\end{proposicion}

\begin{nota} Let us note that property
\begin{equation} \label{eq:gr-ann-ker-phi}
\Gr_{F_T} \ann_{\DX[s]} f^s =\ker \varphi \end{equation} may hold
without assuming any linearity condition. For instance, in examples
\ref{subsec:four-lines}, \ref{subsec:five-lines} the divisor $D$ is
not of differential linear type but property
(\ref{eq:gr-ann-ker-phi}) holds.
\end{nota}

\begin{cuestion} A natural question is whether other locally quasi-homogeneous non-necessarily
 free divisors, for instance arbitrary hyperplane arrangements, are of
 commutative linear type (hence of differential linear type) or not.
 One could try to imitate the proof of the theorem \ref{teo:lqh->clt}, and
 a new question appears: for which locally quasi-homogeneous
 divisors the symmetric algebra of its Jacobian ideal is
 Cohen-Macaulay?
\end{cuestion}

\begin{nota} Let us note that if property (\ref{eq:gr-ann-ker-phi}) holds, let's say at a point $p\in D$, then if $\{P_1,\dots,P_r\}$ is an
involutive basis of $\ann_{\D[s]} f^s$, i.e. $\{\sigma_T(P_1),\dots,\sigma_T(P_r)\}$ is a basis of $\Gr_{F_T} \ann_{\D[s]} f^s = \ker
\varphi_p$, then $\{f,P_1,\dots,P_r\}$ is an involutive basis of $\D[s]f + \ann_{\D[s]} f^s$. This is due to the fact that $\ker \varphi_p$ is a
prime ideal and $(\ker \varphi_p):f = \ker \varphi_p$.
\end{nota}

In what follows we assume that $D$ is a free divisor.
\medskip

The following proposition is proved in \cite{calde_nar_lct_ilc}, prop. (1.6.7). It had been pointed out to us by Torrelli. See also
\cite{simis-alg-free}, cor.~3.12 in the polynomial case.

\begin{proposicion} \label{prop:clt->Kos} If $D$ is a free divisor of commutative linear type
(at $p$), then it is Koszul (at $p$).
\end{proposicion}

 Let $p\in D$ and let us write $\hol=\OXp$ and $\D=\DXp$. Let us consider the following properties:
\begin{enumerate}
\item[\sKp] For some (and hence any) reduced
local equation $f$ of $D$ at $p$ and for some (or any) local basis $\{\delta_1,\dots,\delta_n\}$ at $p$  of the logarithmic vector fields (with
respect to $D$),  one has that
$$ \sigma\left(\delta_1- \frac{\delta_1(f)}{f}s\right),\dots, \sigma\left(\delta_n-
\frac{\delta_n(f)}{f}s\right)$$ is a regular sequence in $\Gr_{F_T}
\D[s] = \Sim_\hol F^1 \D$.
\item[\sK] Property \sKp\ holds at any $p\in D$.
\end{enumerate}

We have the following result \cite{calde_nar_lct_ilc}, prop.
(1.6.2).

\begin{proposicion} \label{prop:thetafs-KF} If $D$ is Koszul   then
$D$ satisfies property  \sK.
\end{proposicion}

\begin{ejemplo} In examples \ref{subsec:four-lines}, \ref{subsec:D-4}, \ref{subsec:D-7}, $D$ satisfies  \sK\ but is not Koszul.
Example \ref{subsec:five-lines} does not satisfy \sK.
\end{ejemplo}

\begin{proposicion} \label{prop:carac-sK} Let us assume that $D$ is an Euler homogeneous free divisor. Let
$f$ be a reduced local equation of $D$ at a point $p$ and
$\{\delta_1,\dots,\delta_n\}$ a local basis at $p$ of the
logarithmic vector fields (with respect to $D$), that we can take in
such a way that $\delta_i(f)=0$ for $i=1,\dots,n-1$ and
$\delta_n(f)=f$. The following properties are equivalent:
\begin{enumerate}
\item[(a)] $D$ satisfies property \sKp\ at $p$.
\item[(b)] $ \sigma(\delta_1),\dots, \sigma(\delta_{n-1})$ is a regular sequence in $\Gr_{F}
\D$.
\end{enumerate}
\end{proposicion}

\begin{prueba} Let  $x_1,\dots,x_n$ be a system of coordinates
centered at $p$ and let us write $\xi_i =\sigma(\partial/\partial
x_i)$ and $\sigma_i = \sigma(\delta_i)$. We have
$$ \Gr_{F_T} \D[s] = \hol[s,\uxi] \supset \Gr_F \D = \hol[\uxi].$$
It is clear by faithful flatness that
$\sigma_1,\dots,\sigma_{n-1},\sigma_n - s$ is a regular sequence in
$\hol[s,\uxi]$ if and only if $\sigma_1,\dots,\sigma_{n-1}$ is a
regular sequence $\hol[\uxi]$.
\end{prueba}

The following proposition (with its corollary) gives a criterion to prove that a free
divisor is pre-Spencer.

\begin{proposicion} \label{prop:crit-sp} Let us assume that $D$ is given by a global
reduced equation $f:X\to\C$ and that the complex $\DX[s] \otimes_{\U
\Thetafs} \SP_{\Thetafs}\OX$ is concentrated in degree $0$, i.e. it
is a $\DX[s]$-resolution of $\DX[s]/\DX[s]\cdot \Thetafs$ (see
(\ref{eq:thetafs})). Then, the following properties are equivalent:
\begin{enumerate}
\item[(a)] The complex $\DX\otimes_{\VO} \SP\OX$ is exact in degree
$-1$.
\item[(b)] The module $\DX[s]/\DX[s]\cdot\Thetafs$ has no
$s$-torsion.
\item[(c)]  The divisor $D$ is pre-Spencer.
\end{enumerate}
\end{proposicion}

\begin{prueba} To prove the proposition we can proceed locally at
any point $p\in D$. Let us write $\hol=\OXp$ and $\D=\DXp$, and let
$\{\delta_1,\dots,\delta_n\}$ be a $\hol$-basis of $\Der(\log D)_p$,
and $\{\delta_1-\alpha_1 s,\dots,\delta_n-\alpha_n s\}$ the
corresponding basis of $\left(\Thetafs\right)_p$, with
$\delta_i(f)=\alpha_i f$. Let us write, for $i<j$,
$[\delta_i,\delta_j] = \sum_{k=0}^n a^{ij}_k \delta_k$ and
$$ R_{ij} =
(a^{ij}_1,\dots,\stackrel{i}{\overbrace{\delta_j+a^{ij}_i}},\dots,
\stackrel{j}{\overbrace{-\delta_i+a^{ij}_j}},\dots,a^{ij}_n)$$ the
corresponding syzygy of $\delta_1,\dots,\delta_n$. Let us also write
$$ {\cal R}_{ij}(s) =
(a^{ij}_1,\dots,\stackrel{i}{\overbrace{\delta_j-\alpha_j s +a^{ij}_i}},\dots, \stackrel{j}{\overbrace{-\delta_i+\alpha_i s
+a^{ij}_j}},\dots,a^{ij}_n)$$ the corresponding syzygy of $\delta_1-\alpha_1 s,\dots,\delta_n-\alpha_n s$. We have ${\cal R}_{ij}(0) = R_{ij}$.
\medskip

\noindent (a) $\Rightarrow$ (b)\quad Let $P(s)\in\D[s]$ such that
$sP(s)\in \D[s]\cdot\Thetafs$:
$$ sP(s) = \sum_{i=1}^n A_i(s)\left( \delta_i-\alpha_i s\right).$$
By taking $s=0$, we obtain $ 0 = \sum_{i=1}^n A_i(0) \delta_i$ and
from (a) we deduce that
$$ (A_1(0),\dots,A_n(0))= \sum_{i<j} c_{ij} R_{ij}.$$
Let us consider
$$ (B_1(s),\dots,B_n(s)) = (A_1(s),\dots,A_n(s)) - \sum_{i<j} c_{ij}
{\cal R}_{ij}(s).$$ By taking $s=0$ again we obtain that all the $B_i(0)$ vanish, so $B_i(s)=s B_i'(s)$,
$$sP(s) = \sum_{i=1}^n A_i(s)\left( \delta_i-\alpha_i s\right) =
\sum_{i=1}^n B_i(s)\left( \delta_i-\alpha_i s\right) = s
\sum_{i=1}^n B_i'(s)\left( \delta_i-\alpha_i s\right)$$and $P(s)\in
\D[s]\cdot\Thetafs$. \medskip

\noindent (b) $\Rightarrow$ (c) \quad Let $C$ be the stalk at $p$ of
$\DX[s] \otimes_{\U \Thetafs} \SP_{\Thetafs}\OX$. It is a free
resolution of $M:=\D[s]/\D[s]\cdot \Thetafs$. The complex $C':= \D
\Lotimes_{\D[s]} M$ can be computed through the $\D[s]$-free
resolution of $\D$ given by $ \D[s] \xrightarrow{s\cdot} \D[s]$ and
so, since $M$ has no $s$-torsion, the complex
$$ C'=\left( \D \Lotimes_{\D[s]} M \right) = \left( \stackrel{-1}{M}
\xrightarrow{s\cdot} \stackrel{0}{M}\right)$$is concentrated in
degree $0$. But $C'$ can be also computed as $C' = \D
\otimes_{\D[s]} C$ and the second complex is nothing but the stalk
at $p$ of $\DX\otimes_{\VO} \SP\OX$, and so it is concentrated in
degree $0$.
\medskip

\noindent (c) $\Rightarrow$ (a)\quad It is trivial.
\end{prueba}

\begin{corolario} \label{cor:crit-SP}
Let us assume that $D$ is given by a global reduced equation
$f:X\to\C$ and that $D$ satisfies property \sK. Then, the following
properties are equivalent:
\begin{enumerate}
\item[(a)] The complex $\DX\otimes_{\VO} \SP\OX$ is exact in degree
$-1$.
\item[(b)] The module $\DX[s]/\DX[s]\cdot\Thetafs$ has no
$s$-torsion.
\item[(c)]  The divisor $D$ is pre-Spencer.
\end{enumerate}
\end{corolario}

\begin{prueba}
It is enough to prove that the property  \sK\ implies that the
complex $L=\DX[s] \otimes_{\U \Thetafs} \SP_{\Thetafs}\OX$ is
concentrated in degree $0$, but this is easily proved by filtering
$L$  in such a way that its graded complex is the Koszul complex
associated with the sequence
$$ \delta_1- \frac{\delta_1(f)}{f}s,\dots, \delta_n-
\frac{\delta_n(f)}{f}s$$ over $\Gr_{F_T} \DX[s] = \Sim_{\OX} F^1 \DX$, where $\{\delta_1,\dots,\delta_n\}$ is a local basis of the logarithmic
vector fields (with respect to $D$). See \cite{calde_nar_lct_ilc}, prop. (1.5.3) for the details.
\end{prueba}

\begin{ejemplo} Corollary \ref{cor:crit-SP} is applied in examples
\ref{subsec:four-lines}, \ref{subsec:D-4}, \ref{subsec:D-7} to
deduce the pre-Spencer property. However, it does not apply to
example \ref{subsec:five-lines}.
\end{ejemplo}

In \cite{torre-overview} the reader can find some related linearity properties and their relationship with the LCT.

\begin{nota} The  property of being pre-Spencer is a non-commutative version
of the notion of regular sequence for $\delta_1,\dots,\delta_n$. However, in the non-commutative case the pre-Spencer property is not inherited
by subsequences, when that makes sense (see example \ref{subsec:five-lines}).
\end{nota}

\subsection{Further linearity conditions: examples and questions}

For a divisor $D\subset X$, the properties of being of commutative linear type or of differential linear type are very restrictive. In dimension
2 we have the following proposition.

\begin{proposicion} \label{prop:carac-2-dim} Let assume that $\dim X =2$, and so $D$ is automatically free. The following properties are equivalent:
\begin{enumerate}
\item[(a)] $D$ is of commutative linear type.
\item[(b)] $D$ is of differential linear type.
\item[(c)] $D$ satisfies the LCT.
\item[(d)] $D$ is Euler homogeneous.
\item[(e)] $D$ is locally quasi-homogeneous.
\end{enumerate}
\end{proposicion}

\begin{prueba} The implication\quad (a) $\Rightarrow$ (b)\quad comes from proposition \ref{prop:clt->dlt}.
The equivalences \quad (c) $\Leftrightarrow$ (d) $\Leftrightarrow$ (e)\quad have been proved in \cite{calde_mon_nar_cas}. The implication \quad
(e) $\Rightarrow$ (a)\quad is a consequence of theorem \ref{teo:lqh->clt}.

For the remaining implication \quad (b) $\Rightarrow$ (d), we use
some results of T. Torrelli (see \cite{calde_nar_lct_ilc}, remark
(1.6.6), c)). We can proceed locally and assume that $D$ has a
reduced equation $f=0$. Since $\dim X=2$, we know that $D$ is free
 and Koszul. Let
$\{\delta_1,\delta_2\}$ be a local basis of the logarithmic
derivations, with $\delta_i(f)=\alpha_i f$. Property (b) means that
$\ann_{\D[s]} f^s$ is generated by $\delta_1-\alpha_1 s, \delta_2
-\alpha_2 s$. But the Bernstein polynomial of $f$ has no integer
roots less than $-1$ (cf. \cite{varchenko_81}), and so $\ann_\D
f^{-1}$ is generated by $\delta_1+\alpha_1, \delta_2 +\alpha_2$. We
can now apply proposition 4.1 and lemma 4.3 of \cite{torre-45-bis}
to conclude that $f$ belongs to the ideal generated by its partial
derivatives, i.e. $D$ is Euler homogeneous.
\end{prueba}

In higher dimension, the relationship between the properties in the above proposition, even for free divisors, is not clear (see conjecture 1.4
in \cite{calde_mon_nar_cas} and \cite{gran-schul-06}). For instance, examples in Section \ref{sec:ex} satisfy the LCT but they are not of
differential linear type. Nevertheless, some generalized linear type conditions appear.
\medskip

Let us denote by $\DX(s):= \C(s)\otimes_\C \DX$ (resp. $\D(s):=\C(s)\otimes_\C \D$), with the filtration $\C(s)\otimes_\C F$, that we will also
denote by $F$.

Let $p\in D$ and let us write $\hol=\OXp$ and $\D=\DXp$. Let us consider the following properties

\begin{enumerate}
\item[\gdltp] For some (or any)\footnote{One easily sees
that this condition does not depend on the choice of the local equation.} reduced local equation $f\in\hol$ of $D$ at $p$, the ideal
$\ann_{\D(s)} f^s$ is generated by order one operators, i.e. (see (\ref{eq:thetafs})) $ \ann_{\D(s)} f^s = \D(s)\cdot \Thetafs$.
\item[\gdlt] Property \gdltp\ holds for any $p\in D$.
\end{enumerate}

It is clear that the property \gdlt\ is an open condition.

\begin{nota}
Since $\ann_{\D[s]} f^s$ is an ideal finitely generated of $\D[s]$, property \gdltp\ holds if and only if there is a non zero polynomial
$\beta(s)\in\C[s]$ such that $\beta(s) \ann_{\D[s]} f^s \subset \ann^{(1)}_{\D[s]} f^s$. A sufficient condition for the existence of a such
polynomial is that the quotient $\ann_{\D[s]} f^s/\ann^{(1)}_{\D[s]} f^s$ is holonomic as $\D$-module. However this condition is not necessary
as shown in example \ref{subsec:four-lines}.
\end{nota}

It is clear that any divisor of differential linear type satisfies
\gdlt.

\begin{ejemplo} The divisors of the examples in \ref{sec:ex} satisfy \gdlt\ but are not of differential linear type.
\end{ejemplo}

Let us assume that $X$ is an open subset of $\C^n$ with coordinates $x_1,\dots,x_n$ and that $D$ has reduced global equation $f:X\to \C$.
Morphism $\varphi$ in (\ref{eq:varphi}) is then given explicitly by
$$\varphi : F(s,\uxi)\in \OX[s,\xi_1,\dots,\xi_n] \mapsto F(f,f'_{\ux})t^d\in \Rees(\OX(f,f'_{x_1},\dots,f'_{x_n}))$$
for each homogeneous polynomial $F$ of degree $d$, where $\xi_i$ denotes the symbol of $\frac{\partial}{\partial x_i}$. Let us denote by
$\ker^{(1)} \varphi$ the ideal of $\OX[s,\uxi]$ generated by $[\ker \varphi]_1\equiv \Thetafs$.
\medskip

To say that $D$ is of commutative linear type means that $\ker \varphi = \ker^{(1)} \varphi$.
\medskip

\numero \label{nume:Delta} By shrinking $X$ if needed, let us take a system of generators of $\Der(\log D)$ on $X$, $ \delta_i = \sum_{j=1}^n
a_{ij} \frac{\partial }{\partial x_j},\quad 1\leq i\leq m$, with $\delta_i(f)=\alpha_i f$. In other words, $(-\alpha_i,a_{i1},\dots,a_{in})$,
$1\leq i\leq m$, is a system of generators of the syzygies of $f,f'_{x_1},\dots,f'_{x_n}$. Clearly the homogeneous polynomials $
\Delta_i=-\alpha_i s + a_{i1} \xi_1 + \cdots +a_{in}\xi_n$, $1\leq i\leq m$, generate $\ker^{(1)} \varphi$.
\medskip

In examples in Section \ref{sec:ex} the divisors are not of differential linear type but they do satisfy the following property: \bigskip

\noindent \gclt\quad\quad There is an integer\quad $N\geq 0$\quad such that\quad $s^N \ker \varphi \subset \ker^{(1)} \varphi$. \bigskip

\noindent In other words, if we consider the map
$$\varphi' : \OX[s,s^{-1},\uxi] \xrightarrow{} \Rees(\OX(f,f'_{\ux}))_{ft}\subset \OX[t]_{ft}$$
induced by $\varphi$, property \gclt\ exactly means that $\ker \varphi'$ is generated by $\ker^{(1)} \varphi$.

Let us note that property \gclt\ does not depend on the particular equation $f$ and so it is a property of $D$. It is clear that any divisor of
commutative linear type satisfies \gclt.

\begin{proposicion} If $D$ satisfies property \gclt, then it is Euler homogeneous.
\end{proposicion}

\begin{prueba} The proof is analogous to the proof of proposition \ref{prop:clt->E-hom}. Let $N>0$ be such that $s^N \ker \varphi \subset \ker^{(1)}
\varphi$. We can proceed locally at each point $p\in D$.

We know that $f$ belongs to the integral closure of the gradient ideal $I=(f'_{x_1},\dots,f'_{x_n})\subset \hol$ (cf. \cite{MR0374482}, \S 0.5,
1)), i.e. there is an integer $d>0$ and elements $a_i \in I^{d-i}$ such that $f^d + a_{d-1} f^{d-1} + \cdots + a_0 =0$. Equivalently, there is a
homogeneous polynomial $F\in \hol[s,\uxi]$ of degree $d>0$ such that $F(f,f'_{x_1},\dots,f'_{x_n})=0$ and $F(s,0,\dots,0)=s^d$. Then $F\in\ker
\varphi_p$ and $s^N F = \sum_{i=1}^m Q_i \Delta_i$, where the $\Delta_i$ have been defined in \ref{nume:Delta}. By taking $\xi_1=\cdots=\xi_n=0$
we deduce that $s^{N+d} = -\sum_i Q_i(s,0\dots,0)  \alpha_i s$ and so at least one of the $\alpha_i$ is a unit, i.e. $f \in
(f'_{x_1},\dots,f'_{x_n})$ and $D$ is Euler homogeneous at $p$.
\end{prueba}

\begin{corolario} If $\dim X=2$ and $D$ satisfies \gclt, then it is of commutative linear type.
\end{corolario}

\begin{nota} Proposition \ref{prop:clt->Kos} suggests that property \gclt\ could imply property \sK, but the divisor in example \ref{subsec:five-lines}
satisfies \gclt\ and not \sK.
\end{nota}

Proposition \ref{prop:clt->dlt} suggests the following question.

\begin{cuestion} We do not know whether the implication \gclt\ $\Rightarrow$\ \gdlt\ is true or not. In the same vein, we do not know whether
there is some relationship between $N$ and the degree of $\beta(s)$ or not, where  $s^N \ker \varphi \subset \ker^{(1)} \varphi$ and
$\beta(s)\ann_{\DX[s]} f^s \subset \ann^{(1)}_{\DX[s]} f^s$ (see examples in \ref{sec:ex})
\end{cuestion}

\begin{cuestion} Concerning property \gdlt, we do not know  any
 ``direct'' algorithm\footnote{Not based on the known algorithms for computing the annihilator of $f^s$ over $\Q[s,\ux,\upartial]$.}
to find generators of the annihilator of $f^s$ over the Weyl algebra $\Q(s)[\ux,\upartial]$, for a given polynomial $f\in\Q[\ux]$.
\end{cuestion}

\begin{cuestion} All the examples of free divisors satisfying the LCT we have been able to compute (see for instance examples in section
\ref{sec:ex}) satisfy properties \gclt\ and \gdlt, but we do not know whether this a general fact or not.
\end{cuestion}

\begin{cuestion} Let $b(s)$ be the Bernstein polynomial of $f$ at a point\footnote{A similar question can be considered in the global polynomial case.} $p\in D$. We have a functional equation
$b(s) f^s = P f^{s+1}$ with $P\in\D[s]$, which means that $b(s) -P f\in \ann_{\D[s]} f^s$. If $D$ satisfies \gdltp\ at $p$, then there is a non
zero polynomial $\beta(s)\in\C[s]$ such that $\beta(s) \ann_{\D[s]} f^s \subset \ann^{(1)}_{\D[s]} f^s$, and so there is a functional equation
$b'(s) f^s = P' f^{s+1}$ with  $b'(s) -P' f\in \ann^{(1)}_{\D[s]} f^s$. In the examples in \ref{sec:ex}, property \gdlt\ is always satisfied (an
also the LCT), but more precise properties hold. Namely, $b(s)$ is a multiple of $\beta(s)$ and the functional equation $b(s) f^s = P f^{s+1}$
for the Bernstein polynomial already satisfies $b(s) -P f\in \ann^{(1)}_{\D[s]} f^s$ without multiplying by $\beta(s)$. We do not know whether
these properties are a consequence of the property \gclt, and eventually of the LCT. On the other hand, a more precise relationship between the
factors of $\beta(s)$ and the tower
$$ \ann_{\D[s]}^{(1)}f^s\varsubsetneq \ann_{\D[s]}^{(2)}f^s \varsubsetneq \cdots\varsubsetneq \ann_{\D[s]}^{(k_0)}f^s=\ann_{\D[s]}f^s$$
appears in examples in section \ref{subsec:fami}.
\end{cuestion}

\begin{cuestion} A natural question is to check whether linear free divisors \cite{gran-mond-nie-schul}, or locally weakly quasi-homogeneous free divisors (LWQH) \cite{cas_gago_hartillo_ucha}
satisfy or not properties \gclt\ or \gdlt. For instance, all the divisors in section \ref{sec:ex} are free LWQH.
\end{cuestion}

\section{Examples}
\label{sec:ex}

In this section we study some examples of free divisors in the affine space $X=\C^n$, for $n=3,4$  with a polynomial equation, from the point of
view of properties \gclt, \gdlt\ and relatives.

As a matter of notation, every time we have a global basis $\{\delta_1,\dots,\delta_n\}$ of the logarithmic derivations, we write
$\sigma_i=\sigma(\delta_i)\in \C[\ux][\uxi] \subset \Gamma(X, \Gr \DX) = \OX(X)[\uxi]$.

The computations are done at the global level of Weyl algebras $\W=\C[\ux,\upartial]$ or $\W[s]$ or polynomial rings, and we have used
\cite{D-mod-M2} and \cite{Plural}.

\subsection{The example of ``four lines'' revisited}
\label{subsec:four-lines}

In this example $X=\C^3$ and $D=\{f=0\}$ with $f=x_1x_2(x_1+x_2)(x_1+x_2x_3)$ (see \cite[Remark 4.2.4]{calde_ens}, \cite[ex.
6.2]{calde_nar_compo} and \cite[4]{calde_mon_nar_cas}). There is a global basis of $\derlogD$, $\{\delta_1,\delta_2,\delta_3\}$, with
$\delta_1(f)=\delta_2(f)=0$, $\delta_3(f)=4f$.
\medskip

Since $(\sigma_1,\sigma_2):\sigma_3 = (\sigma_1,\sigma_2,\tau)$, with $\tau$ a certain homogeneous polynomial  in $\uxi$ of degree $2$, and
$(\sigma_1,\sigma_2):\tau=(x_1,x_2)$, we deduce that $D$ is not Koszul at any point of the line $x_1=x_2=0$, which is strictly contained in the
singular locus of $D$.

The operators $\delta_1, \delta_2, \delta_3 -4s$ form a global $\OX$-basis of $F_T^1\ann_{\DX[s]}(f^s)$ and their symbols with respect to the
total order filtration $\sigma_1,\sigma_2,\sigma_3-4s$ form a regular sequence in $\C[\ux][\uxi,s]$, and so, by flatness, in
$\Gr_{F_T}\D_{X,p}[s]=\OO_{X,p}[\xi_1,\xi_2,\xi_3,s]$ for any $p\in X$. Then, the divisor $D$ satisfies property \sK.
\medskip

The kernel of $ \varphi:\OX[s,\xi_1,\xi_2,\xi_3] \to \Rees(\jac(D))$ is generated by $\sigma_1,\sigma_2,\tau,\sigma_3-4s$.
 Moreover $\left(\ker \varphi\right)\cap \OX[\xi_1,\xi_2,\xi_3] =
\OX[\uxi](\sigma_1,\sigma_2,\tau)$, $\tau \notin
\OX[\uxi,s](\sigma_1,\sigma_2,\sigma_3-4s)$ and $\ker \varphi =
\left(\ker^{(1)} \varphi\right) : f$.
\medskip

\noindent We have $s\tau \in (\sigma_1,\sigma_2,\sigma_3-4s)$ and so
$s \left(\ker \varphi\right) \subset \ker^{(1)}\varphi$. In
particular, $D$ satisfies property \gclt.
\medskip

\noindent On the other hand, it is possible to lift $\tau$ to an operator $T\in \W$ such that $\sigma(T)=\tau$ and $T(f^s)=0$. Then, we deduce
that
$$\ann_{\DX[s]} f^s=
\DX[s](\delta_1, \delta_2,T,\delta_3 -4s),\quad  \ann_{\DX} f^s=
\DX(\delta_1, \delta_2,T),$$ and
$$ \Gr_{F_T}\left(\ann_{\DX[s]} f^s\right) = \ker \varphi, \quad
 \Gr_F \left(\ann_{\DX} f^s\right) = \left(\ker \varphi\right)\cap \OX[\uxi] $$
 (i.e. $D$ satisfies property (\ref{eq:gr-ann-ker-phi}))
and $D$ is not of differential linear type.
\medskip

\noindent We also have $T f \in \W(\delta_1,\delta_2)$ and so $ \ann_{\DX} f^s = \{A\ |\ Af \in \DX(\delta_1,\delta_2) \}$.
\medskip

\noindent We find that $(s+1/2)T \in \W[s](\delta_1, \delta_2,
\delta_3 -4s)$ and so
\begin{equation} \label{eq:4rectas-1}
(s+1/2) \ann_{\DX[s]} f^s \subset \ann^{(1)}_{\DX[s]} f^s,
\end{equation}
i.e. $D$ satisfies property \gdlt. However, the $\DX$-module
$$M=\ann_{\DX[s]}(f^s)/\ann^{(1)}_{\DX[s]}(f^s)$$ is not holonomic: we
have $M = \DX\cdot \overline{P}$ and a syzygy computation  shows
that
$$ \ann_{\DX[s]} \overline{P} = \DX[s] \left(x_1, x_2,
s+1/2\right),\quad \ann_{\DX} \overline{P} = \DX \left(x_1, x_2\right)$$ and so $M$ is not holonomic.
\medskip

\noindent The global Bernstein polynomial of $f$ is
$$b_f(s) = (s+1)^3(s+1/2)(s+3/4)(s+5/4).$$ From (\ref{eq:4rectas-1}) and the fact that the
$b_f(s)$ does not have any integer  root less than $-1$ we deduce that $\ann_{\DX} f^{-1}$ is generated by operators of order 1, even if $D$ is
not of differential linear type.
\medskip

Let us explain now ``why" $D$ is pre-Spencer without testing the definition \ref{def:spencer} itself. Let $p$ be any point of $D$ and let us
write $\OO = \OO_{X,p}$, $\D = \D_{X,p}$.

For the exactness of $\DX\otimes_{\VO} \SP\OX$ in degrees $\neq 0$
we can apply criterion \ref{cor:crit-SP}. We know that
$\sigma_1,\sigma_2,\sigma_3-4s$ is a regular sequence in
$\Gr_{F_T}\D[s]$. We need to prove that
$\D[s]/\ann^{(1)}_{\D[s]}(f^s)$ has no $s$-torsion. For that, let
$Q\in\D[s]$ and suppose that $sQ\in \ann^{(1)}_{\D[s]}(f^s)$. It is
clear that $Q\in \ann_{\D[s]}(f^s)$ and we know from
(\ref{eq:4rectas-1}) that $(s+1/2)Q \in \ann^{(1)}_{\D[s]}(f^s)$,
but
$$sQ,(s+1/2)Q \in \ann^{(1)}_{\D[s]}(f^s) \Rightarrow Q \in
\ann^{(1)}_{\D[s]}(f^s).$$ On the other hand, by specializing (\ref{eq:4rectas-1}) at $s=0$ we deduce that $T\in \D(\delta_1,\delta_2,\delta_3)$
and so
$$ (\sigma_3)+ \left(\ker \Phi_p\right)\cap \OO[\uxi] =
(\sigma_1,\sigma_2,\sigma_3,\tau) \subset \Gr_F
\D(\delta_1,\delta_2,\delta_3).$$ But
$$ \dim \left(\frac{\OO[\uxi]}{(\sigma_1,\sigma_2,\sigma_3,\tau)} \right)
= \dim \left(\frac{\OO[\uxi]}{(\sigma_1,\sigma_2,\tau)} \right) -1 = \dim \Rees(\jac(D)_p)-1 = 3$$ and $\D\otimes_{\D(\log D)} \OO =
\D/\D(\delta_1,\delta_2,\delta_3)$ is holonomic.

Finally, we check that there exists an operator $P$ in the Weyl algebra $\W$ of order $6$ such that the Bernstein functional equation $
b_f(s)f^s= P f^{s+1}$ holds, but for which $b_f(s) - P f$ is not on only in $\ann_{\W[s]} f^s$, but in $\W[s](\delta_1, \delta_2, \delta_3
-4s)$.

\subsection{The example of ``five lines''}
\label{subsec:five-lines}

 In this example $X=\C^4$ and $D=\{f=0\}$
with $f=x_1x_2(x_1+x_2)(x_1+x_2x_3)(x_1+x_2x_4)$. There is a global basis of $\derlogD$,
$\{\delta_1,\delta_2,\delta_3,\delta_4\}$, with\\
$\delta_1(f) = \delta_2(f) = \delta_3(f)=0,\delta_4(f) = 5f$.
\medskip

We have $$\gota:=(\sigma_1,\sigma_2):\sigma_3 \neq
(\sigma_1,\sigma_2),\quad (\sigma_1,\sigma_2):\gota = (x_1,x_2)$$
and so $D$ is neither Koszul nor satisfies property \sKp\ (i.e.
$\sigma_1,\sigma_2,\sigma_3,\sigma_4-5s$ is not a regular sequence
in $\Gr_{F_T}\D_{X,p}[s]=\OO_{X,p}[\uxi,s]$) at any $p$ in the plane
$x_1=x_2=0$, which is strictly contained in the singular locus of
$D$. So we cannot use criterion \ref{cor:crit-SP} to prove that $D$
is pre-Spencer. Nevertheless, a global computation in $\W$ shows
that the complex $\DX\otimes_{\VO} \SP\OX$ is exact in degrees $\neq
0$.
\medskip

The vanishing of $h^{-1} \left(\DX\otimes_{\VO} \SP\OX\right)$
exactly means that the syzygies of
$\delta_1,\delta_2,\delta_3,\delta_4$ are generated by those of
``Spencer type", i.e. those coming from the relations expressing the
brackets $[\delta_i,\delta_j]$, $1\leq i < j\leq 4$, as linear
combinations of the $\delta_k$, $1\leq k\leq 4$ (see the proof of
proposition \ref{prop:crit-sp}). However, the syzygies of
$\delta_1,\delta_2,\delta_3$ are not generated by those of ``Spencer
type". That represents a big difference with regular sequences in
the commutative case.
\medskip

\noindent The kernel of $ \varphi:\OX[s,\uxi] \to \Rees(\jac(D))$ is generated by $\sigma_1,\sigma_2,\sigma_3,\tau_1, \tau_2,
\tau_3,\sigma_4-5s$, where the $\tau_i \in \C[\ux][\uxi]$ are homogeneous in $\uxi$ of degree 2 and $\tau_i \notin \ker^{(1)}\varphi$.
 Moreover $\left(\ker \varphi\right)\cap \OX[\uxi] =
\OX[\uxi](\sigma_1,\sigma_2,\sigma_3,\tau_1,\tau_2,\tau_3)$ and
$\ker \varphi = \left(\ker^{(1)} \varphi\right) : f$.
\medskip

As in example \ref{subsec:four-lines}, we have $s \left(\ker \varphi\right) \subset \ker^{(1)}\varphi$ and so $D$ satisfies property \gclt.
\medskip

\noindent On the other hand, there are operators $ T_i \in \ann_{\DX} f^s$ such that $\sigma(T_i)=\tau_i$, $i=1,2,3$. We deduce that
$\ann_{\DX[s]} f^s= \DX[s](\delta_1, \delta_2,\delta_3, T_1,T_2,T_3,\delta_4 -5s)$, $\ann_{\DX} f^s= \DX(\delta_1, \delta_2, \delta_3,
T_1,T_2,T_3)$ and
$$ \Gr_{F_T}\left(\ann_{\DX[s]} f^s\right) = \ker \varphi, \quad
 \Gr_F \left(\ann_{\DX} f^s\right) = \left(\ker \varphi\right)\cap \OX[\uxi]$$
(i.e. $D$ satisfies property ({\ref{eq:gr-ann-ker-phi})) and $D$ is not of differential linear type.
\medskip

\noindent We also have $T_i f \in \W(\delta_1,\delta_2,\delta_3)$ and so $\ann_{\DX} f^s = \{A\ |\ Af \in \DX(\delta_1,\delta_2,\delta_3) \}$.
\medskip

\noindent We also have the following
\begin{equation} \label{eq:5rectas-1}
(s+2/5) \ann_{\D[s]} f^s \subset \ann^{(1)}_{\D[s]} f^s,
\end{equation}
i.e. $D$ satisfies property \gdlt.
\medskip

\noindent The global Bernstein polynomial of $f$ is
$$b_f(s) = (s+1)^4(s+2/5)(s+3/5)(s+4/5)(s+6/5).$$ From (\ref{eq:5rectas-1}) and the fact
that the $b_f(s)$ does not have any integer  root less than $-1$ we deduce that $\ann_{\DX} f^{-1}$ is generated by operators of order 1, even
if $D$ is not of differential linear type.
\medskip

To conclude that $D$ is Spencer, we have to prove that $\DX\otimes_{\VO} \OX = \DX/\DX(\delta_1,\delta_2,\delta_3,\delta_4)$ is holonomic, but
this can be done in a completely similar way as in the example \ref{subsec:four-lines}.
\medskip

Finally $D$ satisfies the LCT after remark
\ref{nota:carac-LCT-explicit}.

\subsection{The family $D_k=\{f_k=(x_1 x_3+x_2)(x_1^k-x_2^k)=0\}$}
\label{subsec:fami}

This family of examples has been studied in \cite{cas_ucha_exper}.
All of them satisfy the LCT. Here we study their behavior with
respect to properties \gclt, \gdlt\ and relatives.

\subsubsection{The case $k=4$} \label{subsec:D-4}

There is a global basis $\{\delta_1,\delta_2,\delta_3\}$ of $\Der(\log D_4)$ with $\delta_1(f_4)=\delta_2(f_4)=0$ and $\delta_3(f_4)=5 f_4$. We
see that $\sigma_1,\sigma_2$ is a regular sequence and so $D_4$ satisfies property \sK\ (see proposition \ref{prop:carac-sK}) but
$$\gota:=(\sigma_1,\sigma_2):\sigma_3 \neq
(\sigma_1,\sigma_2),\quad (\sigma_1,\sigma_2):\gota = (x_1,x_2)$$ and so $D_4$ is not Koszul at any $p$ in the line $x_1=x_2=0$, which is
strictly contained in the singular locus of $D_4$.
\medskip

\noindent The kernel of $ \varphi:\OX[s,\uxi] \to \Rees(\jac D_4)$ is generated by $\sigma_1,\sigma_2,\sigma_3-5s,\tau$ where $\tau$ is
homogeneous in $\uxi$ of degree $2$ and $\tau\notin \ker^{(1)} \varphi$. Moreover $\left(\ker \varphi\right)\cap \OX[\uxi] =
(\sigma_1,\sigma_2,\tau)$ and $\ker \varphi = \left(\ker^{(1)} \varphi\right) : f_4$.
\medskip

In this example we find that $s^2 \left(\ker \varphi\right) \subset \ker^{(1)}\varphi$ (but $s \left(\ker \varphi\right) \not\subset
\ker^{(1)}\varphi$) and so $D_4$ satisfies property \gclt.
\medskip

In this example, property (\ref{eq:gr-ann-ker-phi}) does not hold, i.e. $\Gr_{F_T} \ann_{\DX[s]} f_4^s \neq \ker \varphi$, or equivalently
$\Gr_F \ann_{\DX} f_4^s \neq \left(\ker \varphi\right)\cap \hol[\uxi]$. \medskip

Let us call $\Gr^{(l)} \ann_{\DX} f_4^s$ the ideal of $\Gr \DX$
generated by $\Gr^i \ann_{\DX} f_4^s$ for $i\leq l$, and
$\ann^{(l)}_{\DX} f_4^s$ the left ideal of $\DX$ generated by $F^l
\ann_{\DX} f_4^s$.
\medskip

We find the following facts:\medskip

\noindent $l=2$: $\Gr^{(2)} \ann_{\DX} f_4^s = (\sigma_1,\sigma_2,
x_1 \tau)$ and $ \ann^{(2)}_{\DX} f_4^s =
\DX(\delta_1,\delta_2,T_1)$ with $\sigma(T_1)=x_1\tau$. We also find
that the commutative relation $s x_1\tau \in \ker^{(1)}\varphi$
lifts to the differential relation $ (s+3/5)T_1 \in
\ann^{(1)}_{\DX[s]} f_4^s$. Moreover,
$$ \ann^{(1)}_{\DX[s]} f_4^s : (s+3/5) =
\DX[s](\delta_1,\delta_2,\delta_3-5s,T_1)= \ann^{(2)}_{\DX[s]} f_4^s.$$

\noindent $l=3$: $\Gr^{(3)} \ann_{\DX} f_4^s = (\sigma_1,\sigma_2,
x_1 \tau,\xi_3\tau)$ and $ \ann^{(3)}_{\DX} f_4^s
=\DX(\delta_1,\delta_2,T_1,T_2)$ with $\sigma (T_2)= \xi_3\tau$. We
have $ (s+2/5)T_2 \in  \ann^{(2)}_{\DX[s]} f_4^s$. Moreover,
$$ \ann^{(2)}_{\DX[s]} f_4^s : (s+2/5) =
 \ann^{(3)}_{\DX[s]} f_4^s.$$
We check that $\ann_{\DX[s]} f_4^s = \ann^{(3)}_{\DX[s]} f_4^s$, $\ann_{\DX} f_4^s = \{A\ |\ Af_4 \in \DX(\delta_1,\delta_2) \}$ and
$\delta_1,\delta_2,T_2$ is a (non involutive) basis of $\ann_{\DX} f_4^s$.

We have then that $ \beta(s) \ann_{\DX[s]} f_4^s \subset \ann^{(1)}_{\DX[s]} f_4^s$, with $\beta(s)=(s+3/5)(s+2/5)$, and $D_4$ satisfies
property \gdlt.
\medskip

We conclude as in example \ref{subsec:four-lines} that $D_4$ is Spencer. \medskip

 The global Bernstein polynomial
of $f_4$ is
$$ b_{f_4}(s)=(s + 1)^3 (s + 6/5)(s + 3/5)(s + 4/5)(s + 2/5)$$
and it is a multiple of $\beta(s)$.

 Although the elements $f,\delta_1,\delta_2,\delta_3-5s$ are not an involutive basis of the ideal of $\DX[s]$ that they generate, it is possible
to check, starting from the commutative relation
$$ s^7 \in \C[\ux][s,\uxi](f,\sigma_1,\sigma_2,\sigma_3-5s),$$
that there exist $P,A_1,A_2,A_3$ in $\C[\ux,\upartial,s]$ such that
$$ b_{f_4}(s) - P f_4 = A_1 \delta_1 + A_2\delta_2+
A_3(\delta_3-5s) \in \ann^{(1)}_{\W[s]} f_4^s,$$ with $P$ of total order equal to $7= \deg b_{f_4}(s)$ (that means that the Bernstein polynomial
satisfies a ``regular" functional equation \cite{kk_houches}).

\subsubsection{The case $k=7$} \label{subsec:D-7}

There is a global basis $\{\delta_1,\delta_2,\delta_3\}$ of $\Der(\log D)$ with $\delta_1(f_7)=\delta_2(f_7)=0$ and $\delta_3(f_7)=8 f_7$.  We
see that $\sigma_1,\sigma_2$ is a regular sequence and so $D_7$ satisfies property \sK\ (see proposition \ref{prop:carac-sK}) but
$$\gota:=(\sigma_1,\sigma_2):\sigma_3 \neq
(\sigma_1,\sigma_2),\quad (\sigma_1,\sigma_2):\gota = (x_1,x_2)$$ and so $D_7$ is not Koszul at any $p$ in the line $x_1=x_2=0$, which is
strictly contained in the singular locus of $D_7$.
\medskip

\noindent The kernel of $ \varphi:\OX[s,\uxi] \to \Rees(\jac D_4)$
is generated by $\sigma_1,\sigma_2,\sigma_3-8s,\tau_1,\tau_2,\tau_3$
where $\tau_i$ is homogeneous in $\uxi$ of degree $i+1$, $1\leq
i\leq 3$, and $\tau_i\notin \ker^{(1)} \varphi$. Moreover
$\left(\ker \varphi\right)\cap \OX[\uxi] =
(\sigma_1,\sigma_2,\tau_1,\tau_2,\tau_3)$ and $\ker \varphi =
\left(\ker^{(1)} \varphi\right) : f_7$.
\medskip

In this example we find that $s^5 \left(\ker \varphi\right) \subset \ker^{(1)}\varphi$ (but $s^4 \left(\ker \varphi\right) \not\subset
\ker^{(1)}\varphi$) and so $D_7$ satisfies property \gclt.
\medskip

In this example, property (\ref{eq:gr-ann-ker-phi}) does not hold either, i.e. $\Gr_{F_T} \ann_{\DX[s]} f_7^s \neq \ker \varphi$, or
equivalently $\Gr_F \ann_{\DX} f_7^s \neq \left(\ker \varphi\right)\cap \hol[\uxi]$. \medskip

\noindent $l=2$: $\Gr^{(2)} \ann_{\DX} f_7^s = (\sigma_1,\sigma_2,
x_1 \tau_1)$ and $ \ann^{(2)}_{\DX} f_4^s =
\DX(\delta_1,\delta_2,T_1)$ with $\sigma(T_1)=x_1\tau_1$. We also
find that the commutative relation $s x_1\tau \in \ker^{(1)}\varphi$
lifts to the differential relation $ (s+3/4)T_1 \in
\ann^{(1)}_{\DX[s]} f_7^s$. Moreover,
$$ \ann^{(1)}_{\DX[s]} f_7^s : (s+3/4) =
\DX[s](\delta_1,\delta_2,\delta_3-8s,T_1)= \ann^{(2)}_{\DX[s]} f_7^s.$$

\noindent $l=3$: $\Gr^{(3)} \ann_{\DX} f_7^s = (\sigma_1,\sigma_2, x_1 \tau_1,\xi_3\tau_1)$ and $ \ann^{(3)}_{\DX} f_7^s
=\DX(\delta_1,\delta_2,T_1,T_2)$ with $\sigma (T_2)= \xi_3\tau_1$. We have $ (s+5/8)T_2 \in  \ann^{(2)}_{\DX[s]} f_7^s$. Moreover,
$$ \ann^{(2)}_{\DX[s]} f_7^s : (s+5/8) =
 \ann^{(3)}_{\DX[s]} f_7^s.$$

\noindent $l=4$: $\Gr^{(4)} \ann_{\DX} f_7^s = (\sigma_1,\sigma_2, x_1 \tau_1,\xi_3\tau_1, x_1 \xi_3 \tau_2)$ and $ \ann^{(4)}_{\DX} f_7^s =
\ann^{(3)}_{\DX} f_7^s + \DX (T_3)$ with $\sigma (T_3)= x_1\xi_3\tau_2$. We have $ (s+1/2)T_3 \in  \ann^{(3)}_{\DX[s]} f_7^s$. Moreover,
$$ \ann^{(3)}_{\DX[s]} f_7^s : (s+1/2) =
 \ann^{(4)}_{\DX[s]} f_7^s.$$

\noindent $l=5$: $\Gr^{(5)} \ann_{\DX} f_7^s = (\sigma_1,\sigma_2, x_1 \tau_1,\xi_3\tau_1, x_1 \xi_3 \tau_2,x_1(x_1\xi_2- 8 \xi_3)\tau_3)$ and
$\ann^{(5)}_{\DX} f_7^s = \DX(\delta_1,\delta_2,T_1,T_2,T_3,T_4)$ with $\sigma (T_4)= x_1(x_1\xi_2- 8 \xi_3)\tau_3$. We have $ (s+3/8)T_4 \in
\ann^{(4)}_{\DX[s]} f_7^s$. Moreover,
$$ \ann^{(4)}_{\DX[s]} f_7^s : (s+3/8) =
 \ann^{(5)}_{\DX[s]} f_7^s.$$

\noindent $l=6$: $\Gr^{(6)} \ann_{\DX} f_7^s = (\sigma_1,\sigma_2, x_1 \tau_1,\xi_3\tau_1, x_1 \xi_3 \tau_2,x_1(x_1\xi_2- 8
\xi_3)\tau_3,(x_1\xi_2- 8 \xi_3)^2\tau_3$ and $ \ann^{(6)}_{\DX} f_7^s =\DX(\delta_1,\delta_2,T_1,T_2,T_3,T_4,T_5)$ with $\sigma (T_5)=
(x_1\xi_2- 8 \xi_3)^2\tau_3$. We have $ (s+1/4)T_5 \in \ann^{(5)}_{\DX[s]} f_7^s$. Moreover,
$$ \ann^{(5)}_{\DX[s]} f_7^s : (s+1/4) =
 \ann^{(6)}_{\DX[s]} f_7^s.$$

We check that $\ann_{\DX[s]} f_7^s = \ann^{(6)}_{\DX[s]} f_7^s$, $\ann_{\DX} f_7^s = \{A\ |\ Af_7 \in \DX(\delta_1,\delta_2) \}$ and
$\delta_1,\delta_2,T_5$ (resp. $\delta_1,\delta_2,T_1,T_2,T_3,T_4,T_5$) is a non involutive (resp. involutive) basis of $\ann_{\DX} f_7^s$.

We have then that $ \beta(s) \ann_{\DX[s]} f_7^s \subset \ann^{(1)}_{\DX[s]} f_7^s$, with $\beta(s)=(s+3/4)(s+5/8)(s+1/2)(s+3/8)(s+1/4)$, and
$D_7$ satisfies property \gdlt.
\medskip

We conclude as in examples \ref{subsec:four-lines} and \ref{subsec:D-4} that $D_4$ is Spencer. \medskip

The global Bernstein polynomial of $f_7$ is
$$ b_{f_7}(s)=
(s+1)^3(s+1/2)(s+7/8)(s+9/8)(s+5/8)(s+3/4)(s+3/8)(s+1/4)$$ and it is a multiple of $\beta(s)$.

In this case we check that $s^{10} \notin (f_7) + \ker \varphi$ and so there is no regular functional equation $b_{f_7}(s)f_7 = P f_7^{s+1}$
with $P\in\W[s]$ of total order $10$.
\medskip

Nevertheless we are able to find a functional equation
\begin{equation} \label{eq:uf}
 b_{f_7}(s) - P f_7 = A_1 \delta_1 + A_2\delta_2+ A_3(\delta_3-8s) \in \ann^{(1)}_{\W[s]} f_7^s,
 \end{equation} with $P$ of total order $15$, but the
computations are quite involved. In fact, we have not been able to find $P$ and the $A_i$ by using directly \cite{D-mod-M2,Plural}. What we have
done is to find a new system of generators of the ideal $I=\W(f_7,\delta_1,\delta_2,\delta_3-8s)$ which ``approximates" an involutive basis by
iterating several times the following process:
\medskip

\noindent (Process) If we start with generators $Q_1,\dots,Q_r$ of $I$, we compute a system of generators of the syzygies of
$\sigma_T(Q_1),\dots,\sigma_T(Q_r)$, and with each of them we check whether their lifting to $\W[s]$ produce a linear combination $Q'$ of the
$Q_i$ whose total symbol is in the ideal $(\sigma_T(Q_1),\dots,\sigma_T(Q_r))$ or not. If the answer is not, we add this $Q'$ to our system of
generators.
\medskip

Once we have constructed a convenient system of generators $Q_1,\dots, Q_s$ containing the original one $f_7,\delta_1,\delta_2,\delta_3-8s$, we
start from the commutative relation $ s^{10} \in \C[\ux][s,\uxi](\sigma_T(Q_1),\dots,\sigma_T(Q_s))$ and we find (\ref{eq:uf}).
\bigskip

Let us note that all examples above are (LWQH) and after remark \ref{nota:LWQH->LCT} we know that they are Spencer and satisfy the LCT.

\bigskip

{\small \noindent Departamento de \'{A}lgebra,
 Facultad de  Matem\'{a}ticas, Universidad de Sevilla, P.O. Box 1160, 41080
 Sevilla, Spain}. \\
{\small {\it E-mail}:  narvaez@algebra.us.es
 }

\end{document}